\documentclass{article}
\newtheorem{pr}{Proposition} 
\newtheorem{lem}{Lemma}
\newtheorem{co}{Corollary} 
\newtheorem{de}{Definition}
 
\newtheorem{theo}{Theorem}

\usepackage{amssymb}
\def\R{\mathbb{R}}
\def\P{\mathbb{P}}
\def\N{\mathbb{N}}
\def\C{\mathbb{C}}
\def\Z{\mathbb{Z}}
\def\Q{\mathbb{Q}}

\def\P{\mathbb{P}}

\def\F{\mathbb{F}}

\begin{document}
\title{From symmetries of the modular tower of genus zero real stable curves to an Euler class for the dyadic circle}
\author{Christophe Kapoudjian\\
        Laboratoire Emile Picard, UMR CNRS 5580\\
        Universit\'e P. Sabatier, Toulouse III\\
        UFR MIG\\
        118, route de Narbonne, 31062 Toulouse Cedex 4\\
        e-mail: ckapoudj@picard.ups-tlse.fr }
\date{February 2001}
\maketitle

\small
\begin{center}
{\bf Abstract}
\end{center}
\indent We build actions of Thompson group $V$ (related to the Cantor set)  and
of
the so-called ``spheromorphism" group of Neretin, on ``towers" of moduli spaces of genus zero
real stable curves. The latter consist of inductive limits of spaces which are
the real parts of the Grothendieck-Knudsen compactification of the usual moduli spaces of punctured
Riemann spheres. By a result of M. Davis, T. Januszkiewicz and
R. Scott, these spaces are aspherical cubical complexes, whose fundamental
groups, the ``pure quasi-braid groups", are some analogues of the classical pure braid groups. By lifting the actions of Thompson
and Neretin groups to the universal covers of the towers, we get new extensions
of both groups by an infinite pure quasi-braid group, and construct what we
call an ``Euler class" for Neretin group, justifying the terminology by
exhibiting an Euler-type cocycle. Further, after introducing the infinite (non-pure) quasi-braid
group, we show that both infinite (non-pure and pure) quasi-braid
groups provide new examples of groups whose classifying spaces, after
plus-construction, are loop spaces.\\ 

      
\normalsize 

The aim of this work is to relate geometrically some discrete groups, such as
Thompson group $V$ (acting on the Cantor set) and the so-called ``spheromorphism``
group of Neretin $N$ (the dyadic analogue of the diffeomorphism group
of the circle, acting on the boundary of the dyadic regular tree ${\cal T}_2$,
cf. \cite{ne}, \cite{kap0}, \cite{kap}, \cite{kap2}), to the moduli spaces
of genus zero curves. More explicitely, let ${\cal M}_{0,n+1}(\C)$ denote the
moduli space $\frac{(\C P^1)^{n+1}\setminus \Delta}{PGL(2,\C)}$, where $\C
P^1$ is the complex projective line, $\Delta$ the thick diagonal. The
Grothendieck-Knudsen-Mumford compactification ${\overline{\cal M}}_{0,n+1}(\C)$ has a concrete
realization as an iterated blow-up along the irreducible components of a
hyperplane arrangement in $\C P^{n-2}$. In \cite{ka1} and \cite{ka2}, the interest of the real
part ${\overline{\cal M}}_{0,n+1}(\R)$ was  revealed, and its topology
has recently been studied in \cite{d-j-s1}, \cite{d-j-s2}, and \cite{de}. The relevance of ${\overline{\cal M}}_{0,n+1}(\R)$ with respect
to Thompson or Neretin groups comes through the common role played by planar trees: ${\overline{\cal M}}_{0,n+1}(\R)$ is a
stratified space, whose strata are labeled by planar trees, while dyadic
planar trees occur in symbols defining elements of Thompson and Neretin groups.\\

In \cite{gr}, P. Greenberg builds a classifying space for Thompson group $F$,
the smallest of Thompson groups, whose
total space is an inductive limit of ``bracelet" spaces ${\cal B}_n$, shown to be
isomorphic, as stratified spaces, to the famous Stasheff associahedra
$K_n$'s. It happens however that this space is too small for supporting
an action of Thompson group $V$. So, the idea
to proceed to an analogous construction for Thompson group $V$ originated in the
observation that a two-sheeted covering space ${\widetilde{\cal
    M}}_{0,n+1}(\R)$ of ${\overline{\cal M}}_{0,n+1}(\R)$ was tiled by $n!$
copies of the associahedron $K_n$, and that its stabilization had the adequate
size of a space which might be dynamically related to $V$.\\

Essentially our construction consists in building two spaces (Proposition 3), ${\widetilde{\cal
    M}}_{0,2^{\infty}}(\R)$ and ${\overline{\cal M}}_{0,3.2^{\infty}}(\R)$ -- which we call the ``towers" of moduli spaces of genus
    zero real curves, the second tower appearing as a real version of a
    stabilized moduli space considered by Kapranov in \cite{ka3} -- ,
defined as inductive limits of spaces ${\widetilde{\cal M}}_{0,n}(\R)$ and
    ${\overline{\cal M}}_{0,n}(\R)$ respectively, and making act cellularly on them not only Thompson group $V$, but also the much
    larger spheromorphism group $N$ (Theorem 1). This result
    may have some interest for two reasons:\\
-- First, $V$ and $N$ are usually
    described as groups almost-acting on trees, or acting on their boundaries,
    which are totally disconnected spaces; though there exist some simplicial
    complexes on which they do act (cf. \cite{br2}, \cite{kap2}), they are
    often ``abstract" topological spaces. Here however, $V$ and $N$ cellularly
    act on a very geometric object (it is an ind-variety of the real algebraic
    geometry).\\
-- Second, $V$ or $N$, as groups of symmetries of the towers of real moduli
    spaces, reveal the $p$-adic (with $p=2$) nature of the latter: indeed, $N$ contains
    the automorphism group $Aut({\cal T}_2)$ of the regular dyadic tree, and in particular the
    group $PGL(2,\Q_2)$, as well as the group of locally analytic
    transformations of the projective line $\Q_2 P^1$, seen as the boundary of
    the dyadic tree (cf. \cite{ne},
    \cite{kap}).\\

One of the towers considered in this text, ${\overline{\cal M}}_{0,3.2^{\infty}}(\R)$, involves the pieces of a {\it cyclic} operad
    $\{{\cal M}(n)=\overline{\cal M}_{0,n+1}(\R),\; n\geq 1\}$, in the sense
    of Ginzburg-Kapranov, cf. \cite{gi-ka}. Now between
    Thompson group $F$, acting on the interval $[0,1]$, Thompson group $V$,
    acting on the Cantor set, stands Thompson group $T$, related to
    cyclicity through its realization as a homeomorphism group of the circle
    (cf. \cite{gh-se}). We show that $T$ acting on the ``cyclic" tower ${\overline{\cal
    M}}_{0,3.2^{\infty}}(\R)$ stabilizes a fundamental tile, isomorphic to an
    infinite associahedron $K_{\infty}$, making the other tiles turn around
    $K_{\infty}$. We prove (Theorem 2) that the stabilizer of $K_{\infty}$ under the
    action of $N$ is the
    semi-direct product $Out(T)\triangleright T$, where $Out(T)$, the exterior automorphism
    group of $T$, is known to be the cyclic order 2 group. It is remarkable that passsing from one tile ($K_{\infty}$)
    to the whole tower ${\overline{\cal M}}_{0,3.2^{\infty}}(\R)$, the size of
    the symmetry group, from countable ($T$), becomes uncountable ($N$).\\

 Further, by lifting the actions to
the universal covering spaces, we deduce the existence of non-trivial extensions of
$V$ and $N$: by an infinite ``pure quasi-braid group" -- denoted $PJ_{2^{\infty}}$
    through this text -- , an
inductive limit of the fundamental groups $PJ_n$ of the spaces ${\widetilde{\cal M}}_{0,n+1}(\R)$, for the first action, and by a
group $Q_{3.2^{\infty}}$ for the second one, an inductive limit of the
fundamental groups $Q_n$ of the spaces ${\overline{\cal M}}_{0,n}(\R)$
(In fact, $PJ_{2^{\infty}}$ and $Q_{3.2^{\infty}}$ happen to be isomorphic). In their
definitions, though not in essence, they are somewhat analogous to the pure braid groups $P_n$, this
explains the terminology (\cite{d-j-s2}, \cite{de}). In particular, the resulting extension:
$$(*)\;\;1\rightarrow PJ_{2^{\infty}} \longrightarrow \tilde{\cal A}_N \longrightarrow
N\rightarrow 1 $$
leads to a new non-trivial central extension of $N$
with kernel $\Z/2\Z$ (Theorem 3), after a certain abelianization process of
the kernel. Viewing the boundary $\partial{\cal
  T}_2\cong \Q_2 P^1$ as the ``dyadic circle", we suggest thinking of the resulting class as the
analogue for $N$ of the Euler class of the homeomorphism group of the circle
$Homeo^+(S^1)$, or of Thompson group $T$ (cf. \cite{gh-se}), and
strengthen the analogy by exhibiting an Euler-type cocycle related to the
central extension (Theorem 4).\\  

The first section is devoted to the description of the real moduli spaces
${\overline{\cal M}}_{0,n+1}(\R)$, especially the combinatorics of their
stratifications. Section 2 introduces the infinite moduli space
${\widetilde{\cal M}}_{0,2^{\infty}}(\R)$ and ${\overline{\cal
    M}}_{0,3.2^{\infty}}(\R)$ (the ``towers", cf. Proposition 3). After proving
that Neretin group $N$ (and so Thompson group $V$)  acts upon both towers (Theorem 1) and
computing the stabilizer of a maximal tile (Theorem 2), we define the resulting extensions
 by the infinite
pure quasi-braid groups $PJ_{2^{\infty}}$ and $Q_{3.2^{\infty}}$. In section
3, the infinite (non-pure) quasi-braid group
$J_{2^{\infty}}$ is introduced (Proposition 4), coming in a
short exact sequence (a restriction of the extension $(*)$ above):
$$1\rightarrow PJ_{2^{\infty}} \longrightarrow J_{2^{\infty}}
\longrightarrow \Sigma_{2^{\infty}}\rightarrow 1,$$
where $\Sigma_{2^{\infty}}$ is an inductive limit of
symmetric groups. A ``stable" length is defined on $J_{2^{\infty}}$, related to
a modified word metric on the quasi-braid groups, and is the key for building the Euler class of
$N$ (Theorem 3) and defining the associated 2-cocycle (Theorem 4). In section
4 we claim that the Euler class of $N$ restricted to $PGL(2,\Q_2)$ is
trivial (Theorem 5); we show finally that the classifying spaces of the
infinite quasi-braid groups are homologically equivalent to loop spaces.\\

\noindent{\bf Acknowledgements}: This work is certainly the fruit of a research group session whose motivation
was to understand the relations between Thompson groups and geometric
operads. The author is grateful to V. Sergiescu, as the initiator of this
research group first, and for many discussions concerning the present paper
also. He is thankful to T. Januszkiewicz, Yu. Neretin, R. Scott, and C. Roger
for commentaries or helpful remarks concerning the first version of the paper.   

\section{The compactified real moduli space ${\overline{\cal M}}_{0,n+1}(\R)$}

\subsection{The space of real $(n+1)$-stable curves ${\overline{\cal M}}_{0,n+1}(\R)$}

\begin{de}
Let $k$ be the field $\C$ or $\R$. A genus zero $(n+1)$-stable curve is an algebraic
curve ${\cal C}$ on the field $k$ with $(n+1)$
punctures $x_0,\ldots,x_n$ such that
\begin{enumerate}
\item each irreducible component of ${\cal C}$ is isomorphic to a projective line
$\P_k^1$, and each double point of ${\cal C}$ is ordinary;
\item the graph of ${\cal C}$ is a tree;
\item each component of ${\cal C}$ has at least three points, double or marked.
\end{enumerate}
\end{de}

The graph of a curve $({\cal C},x_0,\ldots,x_n)$ is defined as follows: The
leaves 
(or 1-valent vertices) $A_0,\ldots,A_n$ are in correspondence with the marked points of the curve,
$x_0,\ldots,x_n$, whereas the internal vertices $v_1,\ldots,v_k$ are in bijection with the irreducible
components ${\cal C}_1,\ldots,{\cal C}_k $. There is an internal edge
$[v_iv_j]$ if ${\cal C}_i$ and ${\cal C}_j$ contain the same double
point. Terminal edges $[A_i v_{j_i}]$ correspond to pairs $(x_i,{\cal C}_j)$ with
$x_i\in {\cal C}_{j_i}$.\\
From now on, we shall follow the convention of \cite{gi-ka}: edges will be
oriented in such a way that $[A_0 v_{j_0}]$ is the output edge, and $[A_i
v_{j_i}]$, $i=1,\ldots,n$ are the input edges.   

\subsection{Explicit construction of ${\overline{\cal M}}_{0,n+1}(\R)$}  
 We find it useful to recall the construction of ${\overline{\cal
 M}}_{0,n+1}(\R)$, using mixed sources, namely \cite{be-gi}, \cite{gi-ka},
 \cite{ka1}, \cite{ka2} and \cite{de}:\\  

(1) ${\cal M}_{0,n+1}(\R)=\frac{(\R P^1)^{n+1}\setminus
    \Delta}{PGL(2,\R)}\cong \frac{\R^n\setminus
    \tilde{\Delta}}{Aff(2,\R)}\hookrightarrow \P ^{n-2}_{\R}:=
\P\{(a_1,\ldots,a_n)\in
    \R^n:\;\sum_{i=1}^n a_i=0\}$, where $\Delta$ (resp. $\tilde{\Delta}$) is the
    thick diagonal of $(\R P^1)^{n+1}$ (resp. $\R^n$), and $Aff(2,\R)$ is the
    affine group acting on $\R^2$. The first isomorphism is obvious; as for
    the embedding in the $(n-2)$-dimensional projective space, it is induced
    by the map sending $(x_1,\ldots,x_n)\in\R^n\setminus\tilde{\Delta}$ to
    $(a_1,\ldots,a_n)$, with $a_i=x_i-\frac{1}{n}\sum_{j=1}^n x_j$.\\
The image of ${\cal M}_{0,n+1}(\R)$ in $\P^{n-2}_{\R}$ is the complementary of the
    union of hyperplanes $H_{i,j}:\, a_i=a_j$. This hyperplane arrangement is
    called the {\it braid arrangement}. So ${\cal M}_{0,n+1}(\R)$
    is diconnected into $\frac{n!}{2 }$ connected components \footnote{Note the
    difference with the complex case: ${\cal M}_{0,n+1}(\C)$ is connected}, which are the
    projective Weyl chambers $W^{\sigma}: a_{\sigma(1)}<a_{\sigma(2)}<\ldots<
    a_{\sigma(n)}$, where $\sigma$ belongs to the symmetric group
    $\Sigma_n$. Projectivisation introduces an identification of the chambers
    $W^{\sigma}$ and $W^{\sigma\omega}$, where $\omega$ is the permutation $\left(\begin{array}{cccc}
1 & 2 & \ldots & n\\
n & n-1 & \ldots & 1
\end{array}\right).$\\

(2) Denote by ${\mathfrak B}_0$ the set of $n$  points $p_i$: 
    $a_1=\ldots=\widehat{a_i}=\ldots=a_n$, ${\mathfrak B}_1$ the set of lines $(p_ip_j)$:
    $a_1=\ldots\hat{a_i}=\ldots=\hat{a_j}=\ldots=a_n$, and more generally, by
    ${\mathfrak B}_k$ the set of $k$-planes
    $a_1=\ldots\widehat{a_{i_1}}=\ldots=\widehat{a_{i_{k+1}}}=\ldots=a_n$, 
$i_1<\ldots<i_{k+1}$, for $k=0,
    \ldots,n-3$ (they are the irreducible components, in the sense of Coxeter
    groups: this will become clear in \S1.3) . Along  components of ${\mathfrak B}_k$, hyperplanes
    $H_{i,j}$ do not meet transversely. So, the construction is the following:
    points of ${\mathfrak B}_0$ are first
    blown-up, and we get the blow-down map $X_1\stackrel{\pi_1}{\longrightarrow}\P^{n-2}_{\R}$. The proper transforms
    of the lines $(p_ip_j)$ (i.e. the closures of
    $\pi_1^{-1}((p_ip_j)\setminus\{p_i,p_j\})$) become transverse in $X_1$,
    consequently they can be blown-up in
    any order to produce
    $X_2\stackrel{\pi_2}{\longrightarrow}X_1\stackrel{\pi_1}{\longrightarrow}\P^{n-2}_{\R}$;
    again the proper transforms of the planes $(p_i p_j p_k)$ become
    transverse in $X_2$, and 
    are blown-up in any order. Finally we get the composition of blow-down maps
$${\overline{\cal
    M}}_{0,n+1}(\R):=X_{n-2}\stackrel{\pi_{n-2}}{\longrightarrow}\ldots\stackrel{\pi_3}{\longrightarrow}
X_2\stackrel{\pi_2}{\longrightarrow} X_1\stackrel{\pi_1}{\longrightarrow}\P^{n-2}_{\R},$$
where the last blow-up is useless, since blowing-up along hypersurfaces
does not change the manifold. We sall denote by $p:{\overline{\cal
    M}}_{0,n+1}(\R)\longrightarrow \P^{n-2}_{\R}$ the composition of the iterated blow-ups.\\

(3) Each blow-up along a smooth algebraic submanifold produces a smooth exceptional
    divisor in the new manifold, which is isomorphic to the projective normal
    bundle over the submanifold. We denote by $\widehat{{\mathfrak B}}_k$ the
    set of
    proper transforms in ${\overline{\cal M}}_{0,n+1}(\R)$ of the exceptional
    divisors produced by blowing-up the (proper transforms in $X_k$ of the)
    components of ${\mathfrak B}_k$: they are smooth irreducible
    hypersurfaces of ${\overline{\cal M}}_{0,n+1}(\R)$, and meet
    transversely. So the union $\cup_{k=0}^{n-3} \widehat{{\mathfrak B}}_k$ is the set of
    irreducible components of a {\it normal crossing
    divisor} $\hat{D}$.\\

(4) The real algebraic variety ${\overline{\cal M}}_{0,n+1}(\R)$ is stratified
    in the following way: the closures of the codimension $k$ strata -- or
    closed strata -- are the
    non-empty intersections of $k$ irreducible components of the divisor $\hat{D}$,
    with one of the $\frac{n!}{2}$ closures of the preimages
    $p^{-1}(W^{\sigma})$.  A stratum -- or open
    stratum, though it is not an open set -- is then defined as the complementary of a all closed
    strict substrata in a given closed stratum. With this terminology, an open
    (resp. closed) stratum is a smooth manifold (resp. with boundary) and an
    open (resp. closed) cell.\\
The construction of the complex algebraic variety ${\overline{\cal
    M}}_{0,n+1}(\C)$ is similar, but the complex codimension $k$ closed strata
    are simply the non-empty intersections of $k$ irreducible components of
    the divisor $\hat{D}$, and because the divisor $\hat{D}$
    is normal, they are indeed smooth closed manifolds (without boundary), but
    not cells.\\

(5) Each codimension $k$ stratum ${\cal M}$ is coded by a {\it planar tree}
that we now define:\\
If $\overline{\cal M}=\overline{p^{-1}(W^{\sigma})}\cap D_1\cap\ldots\cap D_k$
is non-empty, then all the components $D_{\alpha}$, $\alpha=1,\ldots k$ are
produced by a blow-up along a set which {\it must be} of the form
$a_{\sigma(i)}=a_{\sigma(i+1)}=\ldots=a_{\sigma(j)}$, with $1\leq i<j\leq n$.\\
Let $T(\sigma, (i,j))$ denote the planar rooted tree with one output edge
and $n$ input edges, labeled from the left to the right by
$\sigma(1),\ldots,\sigma(n)$, with a unique internal edge:

$$
\begin{picture}(0,0)%
\includegraphics{strat.pstex}%
\end{picture}%
\setlength{\unitlength}{3108sp}%
\begingroup\makeatletter\ifx\SetFigFont\undefined%
\gdef\SetFigFont#1#2#3#4#5{%
  \reset@font\fontsize{#1}{#2pt}%
  \fontfamily{#3}\fontseries{#4}\fontshape{#5}%
  \selectfont}%
\fi\endgroup%
\begin{picture}(3072,4389)(3241,-4165)
\put(4726,-4111){\makebox(0,0)[lb]{\smash{\SetFigFont{9}{10.8}{\rmdefault}{\mddefault}{\updefault}figure 1}}}
\put(5401,-3211){\makebox(0,0)[lb]{\smash{\SetFigFont{9}{10.8}{\rmdefault}{\mddefault}{\updefault}. . .}}}
\put(5401,-1861){\makebox(0,0)[lb]{\smash{\SetFigFont{9}{10.8}{\rmdefault}{\mddefault}{\updefault}. . .}}}
\put(4141,-1861){\makebox(0,0)[lb]{\smash{\SetFigFont{9}{10.8}{\rmdefault}{\mddefault}{\updefault}. . .}}}
\put(3556,-2086){\makebox(0,0)[lb]{\smash{\SetFigFont{9}{10.8}{\rmdefault}{\mddefault}{\updefault}$A_{\sigma(1)}$}}}
\put(6211,-2086){\makebox(0,0)[lb]{\smash{\SetFigFont{9}{10.8}{\rmdefault}{\mddefault}{\updefault}$A_{\sigma(n)}$}}}
\put(3241,-3391){\makebox(0,0)[lb]{\smash{\SetFigFont{9}{10.8}{\rmdefault}{\mddefault}{\updefault}$A_{\sigma(i)}$}}}
\put(6031,-3391){\makebox(0,0)[lb]{\smash{\SetFigFont{9}{10.8}{\rmdefault}{\mddefault}{\updefault}$A_{\sigma(j)}$}}}
\put(4411,-3391){\makebox(0,0)[lb]{\smash{\SetFigFont{9}{10.8}{\rmdefault}{\mddefault}{\updefault}$A_{
\sigma(i+1)}$}}}
\put(4861, 29){\makebox(0,0)[lb]{\smash{\SetFigFont{9}{10.8}{\rmdefault}{\mddefault}{\updefault}$A_0$}}}
\end{picture}
$$

Define a {\it contraction} of a tree as the operation consisting in collapsing an internal edge on a single
vertex. Introduce the partial order on the set of $n$-planar trees: $T\leq T'$
if $T'$ is obtained from $T$ by a sequence of contractions.

\begin{de}
The tree attached to the closed stratum $\overline{\cal M}=\overline{p^{-1}(W^{\sigma})}\cap
D_1\cap\ldots\cap D_k$ (or to the stratum ${\cal M}$) is defined to be
$$T(\sigma)=inf_{\alpha=1,\ldots,k}{T(\sigma,(i_{\alpha},j_{\alpha}))}.$$
\end{de}

{\bf Notation:} We shall use two different notations, not to be confused:
 $T(\sigma)$, with $\sigma\in \Sigma_n$, will denote an $n$-planar tree with
 leaves labeled from $\sigma(1)$ to $\sigma(n)$, from the left to the right,
 whereas $(T,\sigma)$ will refer to the same tree, but with the canonical
 labeling of the leaves, from $1$ (on the left) to $n$ (on the right), coupled with the same permutation
 $\sigma$.\\

From the preceding description, the following facts are obvious:\\

{\bf Fact 1:} The stratum $\overline{\cal M}=\overline{p^{-1}(W^{\sigma})}\cap D_1\cap\ldots\cap D_k$
is non-empty if and only if the collection of sets
$S_{\alpha}=\{i_{\alpha},i_{\alpha}+1,\ldots,j_{\alpha}\}$, attached to each
$D_{\alpha}$ as explained above, for $\alpha=1,\ldots,k$, is {\it nested} in
the following sense:\\
\begin{center}
$\forall \alpha,\;\beta,$ either $S_{\alpha}\cap S_{\beta}=\emptyset$,
$S_{\alpha}\subset S_{\beta}$ or $S_{\beta}\subset S_{\alpha}$.
\end{center}

{\bf Fact 2:} Denoting the stratum by ${\cal M}(T,\sigma)$,
we observe that $codim\, {\cal M}(T,\sigma) = card\{internal\; edges\}= card\; Vert(T)\setminus\{root\}$. Indeed, the set of internal edges of $T$ is in
one-to-one correspondence with the codimension $1$ components $D_1,\ldots,D_k$
containing
$\overline{{\cal M}(T,\sigma)}$.\\

Remark: In view of \S 1.1, it is enlighting to understand each stratum ${\cal M}(T,\sigma)$ as the
   set of stable curves whose associated planar tree is $T(\sigma)$. Then Fact
   1 becomes completely clear.

\subsection{Coxeter group formulation, following Davis, Januszkie\-wicz and Scott} 

Following \cite{d-j-s2}, we now formulate the condition guaranteeing a
collection is nested in group theoretic terms, namely in terms of the symmetric group $\Sigma_n$.\\
Denote by $S=\{\sigma_1,\ldots,\sigma_{n-1}\}$ the set of canonical Coxeter generators of
 $\Sigma_n$: $\sigma_i$ is the transposition $(i,i+1)$. Let $(T,\sigma)$ be a
 labeled planar $n$-tree: to each vertex $v$ except the root of $T$, a
proper subset
$T_v=\{\sigma_i,\sigma_{i+1},\ldots,\sigma_{j-1}\}$ of $S$ is associated, corresponding to a
connected subgraph $G_{T_v}$ of the Coxeter graph of $\Sigma_n$: $i,i+1,\ldots,j$ are
the labels of the leaves descending from $v$. The collection ${\cal
  T}=\{T_v,\, v\in Vert(T)\setminus\{root\}\}$ is a {\it nested collection}
in the following sense:

\begin{de}[\cite{d-j-s2}]
A collection ${\cal T}$ of proper subsets of $S$ will be called nested if the
following conditions hold:
\begin{enumerate}
\item The Coxeter subgraph $G_T$ is connected for all $T\in{\cal T}$.
\item For any $T,T'\in {\cal T}$, either $T\subset T'$, $T'\subset T$, or
  $G_{T\cup T'}$ is not connected.
\end{enumerate}
\end{de}

It is now clear that there is a bijection $T\leftrightarrow {\cal T}$ between the set of planar $n$-trees  and the set of nested
  collections of the symmetric group $\Sigma_n$. We shall later use the induced
  correspondence $(T,\sigma)\leftrightarrow ({\cal T},\sigma)$.\\

\subsection{Combinatorics of the stratification of a $\overline{\cal
  M}_{0,n+1}(\R)$ and of its two-sheeted cover $\widetilde{\cal M}_{0,n+1}(\R)$}

{\bf 1.4.1} Each $K^{\sigma}=\overline{p^{-1}(W^{\sigma})}$ is combinatorially
    isomorphic to the Stasheff associahedron $K_n$ (cf. \cite{st}). This is so because
    the stratifications of both objects are the same. Of course, strata
    labeled by different planar trees may coincide in the space
    ${\overline{\cal M}}_{0,n+1}(\R)$. In \cite{de}, the rules of
    identification of faces are expressed
    in the language of polygons, which is well-adapted because the operad
    $\{{\cal M}(n)={\overline{\cal M}}_{0,n+1}(\R),\; n>1\}$ is cyclic in the
    sense of Getzler-Kapranov (cf. \cite{ge-ka}, the action of $\Sigma_n$ on
    ${\cal M}(n)$ extends in an obvious way to $\Sigma_{n+1}$). However, we prefer
    translating them in terms of trees, because of their relevance with
    respect to Thompson groups. Moreover, instead of ${\overline{\cal
    M}}_{0,n+1}(\R)$, we shall consider a two-sheeted covering ${\widetilde{\cal M}}_{0,n+1}(\R)$, that we now define:\\

Consider the two-sheeted covering $S^{n-2}\longrightarrow
  \P_{\R}^{n-2}$, where $S^{n-2}$ is the $(n-2)$-dimensional unit sphere of
  $\{(a_1,\ldots,a_n)\in\R^n:\;\sum_{i=1}^n a_i=0\}$.\\
Apply now the process of iterated blow-ups described in \S 1.1 to
  $S^{n-2}$ with its lifted braid arrangement, and denote by ${\widetilde{\cal
      M}}_{0,n+1}(\R)$ the resulting space. This yields a commutative diagram

$${\widetilde{\cal
    M}}_{0,n+1}(\R)\longrightarrow\overline {\cal
  M}_{0,n+1}(\R)$$

$$S^{n-2}\longrightarrow  \P_{\R}^{n-2}$$

where the horizontal arrows are the obvious two-sheeted covering maps, and the vertical
ones are the blow-down maps. Since $S^{n-2}$ is tiled by $n!$ Weyl
chambers, the covering ${\widetilde{\cal M}}_{0,n+1}(\R)$ will be tiled by
$n!$ copies of the associahedron $K_n$.\\

\subsubsection*{1.4.2 Combinatorics of the stratification of ${\widetilde{\cal
      M}}_{0,n+1}(\R)$.} Let $T(\sigma)$ be an $n$-planar tree, and $v$ an internal vertex other than
the root. The tree disconnects in $v$
into 2 subtrees, the ``upper" one (with the same output edge) and the ``lower" one
(with root $v$). Proceed now to a reflection of the lower tree -- so that the
left-to-right labeling of its input edges is inversed -- before
glueing both pieces back together to form a new planar $n$-tree, denoted
${\tilde{\nabla}}_v T(\sigma)$. We shall also denote by ${\tilde{\nabla}}_v
(T,\sigma)$ the couple $(T',\sigma')$ such that
$T'(\sigma')={\tilde{\nabla}}_v T(\sigma)$ (cf. figure 2 for an example). 

\begin{pr}[following \cite{de}]
Two strata ${\cal M}(T,\sigma)$ and ${\cal M}(T',\sigma')$ in the double cover
${\widetilde{\cal M}}_{0,n+1}(\R)$ coincide if and only if
$(T,\sigma)\sim (T',\sigma')$, where $\sim $ is the equivalence relation
generated by $(T,\sigma)\sim{\tilde{\nabla}}_v(T,\sigma)$ for any vertex
distinct from the root. In other words, $(T,\sigma)\sim (T',\sigma')$ if and
only if $\exists v_1\in Vert(T)$, $\exists v_i\in
Vert({\tilde{\nabla}}_{v_{i-1}}\ldots{\tilde{\nabla}}_{v_1} T(\sigma))$, $i=2,\ldots,r$,
such that $(T',\sigma')={\tilde{\nabla}}_{v_r}\ldots{\tilde{\nabla}}_{v_1} (T,\sigma)$. 
\end{pr}

Remark: It is easy to see that it can be supposed that $v_1,\ldots,v_r$
are vertices of the same initial tree $T$, and if $v_i$ is a descendent of
$v_j$, then $i\leq j$.\\

Proof: The rule $(T,\sigma)\sim{\tilde{\nabla}}_v(T,\sigma)$ comes from the
projectivisation of the normal bundle over the blown-up components. Indeed, suppose we
blow-up the component $a_{\sigma(i)}=a_{\sigma(i+1)}=\ldots=a_{\sigma(j)}$;
the equations of the blow-up are:
$$\lambda_k(a_{\sigma(l)}-a_{\sigma(i)})=\lambda_l(a_{\sigma(k)}-a_{\sigma(i)}),\;\;k,l=i+1,\ldots,j$$
$$[\lambda_{i+1}:\ldots:\lambda_{j}]\in\R P^{j-i-1},\;\;(a_1,\ldots,a_n)\in
S_{\R}^{n-2},$$
and because $[\lambda_{i+1}:\ldots:\lambda_{j}]=[-\lambda_{i+1}:\ldots:-\lambda_{j}]$,
it follows that a point of the exceptional divisor close to the cell
$$a_{\sigma(1)}<a_{\sigma(2)}<\ldots<a_{\sigma(i)}<\ldots<a_{\sigma(j)}<\ldots<a_{\sigma(n)}$$
is close also to the cell  
$$a_{\sigma(1)}<a_{\sigma(2)}<\ldots<a_{\sigma(j)}<\ldots<a_{\sigma(i)}<\ldots<a_{\sigma(n)}.$$
This gives the identification rules for codimension 1 strata, and since the other strata are intersections of codimension 1 strata, the complete
rule follows. $\;\;\square$\\
 
$$
\begin{picture}(0,0)%
\includegraphics{recol.pstex}%
\end{picture}%
\setlength{\unitlength}{3108sp}%
\begingroup\makeatletter\ifx\SetFigFont\undefined%
\gdef\SetFigFont#1#2#3#4#5{%
  \reset@font\fontsize{#1}{#2pt}%
  \fontfamily{#3}\fontseries{#4}\fontshape{#5}%
  \selectfont}%
\fi\endgroup%
\begin{picture}(6885,3039)(316,-2815)
\put(1126,-1186){\makebox(0,0)[lb]{\smash{\SetFigFont{9}{10.8}{\rmdefault}{\mddefault}{\updefault} $v$}}}
\put(4996,-1186){\makebox(0,0)[lb]{\smash{\SetFigFont{9}{10.8}{\rmdefault}{\mddefault}{\updefault}$v$}}}
\put(1891,-646){\makebox(0,0)[lb]{\smash{\SetFigFont{9}{10.8}{\rmdefault}{\mddefault}{\updefault}$a_1$}}}
\put(5941,-646){\makebox(0,0)[lb]{\smash{\SetFigFont{9}{10.8}{\rmdefault}{\mddefault}{\updefault}$a_1$}}}
\put(766,-1276){\makebox(0,0)[lb]{\smash{\SetFigFont{9}{10.8}{\rmdefault}{\mddefault}{\updefault}$A_1$}}}
\put(4636,-1366){\makebox(0,0)[lb]{\smash{\SetFigFont{9}{10.8}{\rmdefault}{\mddefault}{\updefault}$A_1$}}}
\put(3376,-2761){\makebox(0,0)[lb]{\smash{\SetFigFont{9}{10.8}{\rmdefault}{\mddefault}{\updefault}figure 2}}}
\put(5041, 29){\makebox(0,0)[lb]{\smash{\SetFigFont{9}{10.8}{\rmdefault}{\mddefault}{\updefault}$A_0$}}}
\put(1216, 29){\makebox(0,0)[lb]{\smash{\SetFigFont{9}{10.8}{\rmdefault}{\mddefault}{\updefault}$A_0$}}}
\put(3736,-556){\makebox(0,0)[lb]{\smash{\SetFigFont{9}{10.8}{\rmdefault}{\mddefault}{\updefault}${\tilde{\nabla}}_v$}}}
\put(2116,-1816){\makebox(0,0)[lb]{\smash{\SetFigFont{9}{10.8}{\rmdefault}{\mddefault}{\updefault}$a_3$}}}
\put(2746,-1636){\makebox(0,0)[lb]{\smash{\SetFigFont{9}{10.8}{\rmdefault}{\mddefault}{\updefault}$a_4$}}}
\put(2836,-1276){\makebox(0,0)[lb]{\smash{\SetFigFont{9}{10.8}{\rmdefault}{\mddefault}{\updefault}$a_5$}}}
\put(2881,-646){\makebox(0,0)[lb]{\smash{\SetFigFont{9}{10.8}{\rmdefault}{\mddefault}{\updefault}$a_6$}}}
\put(6976,-646){\makebox(0,0)[lb]{\smash{\SetFigFont{9}{10.8}{\rmdefault}{\mddefault}{\updefault}$a_6$}}}
\put(6166,-1411){\makebox(0,0)[lb]{\smash{\SetFigFont{9}{10.8}{\rmdefault}{\mddefault}{\updefault}$a_5$}}}
\put(316,-2671){\makebox(0,0)[lb]{\smash{\SetFigFont{9}{10.8}{\rmdefault}{\mddefault}{\updefault}$A_2$}}}
\put(946,-2671){\makebox(0,0)[lb]{\smash{\SetFigFont{9}{10.8}{\rmdefault}{\mddefault}{\updefault}$A_3$}}}
\put(1216,-2086){\makebox(0,0)[lb]{\smash{\SetFigFont{9}{10.8}{\rmdefault}{\mddefault}{\updefault}$A_4$}}}
\put(1801,-2086){\makebox(0,0)[lb]{\smash{\SetFigFont{9}{10.8}{\rmdefault}{\mddefault}{\updefault}$A_5$}}}
\put(1711,-1276){\makebox(0,0)[lb]{\smash{\SetFigFont{9}{10.8}{\rmdefault}{\mddefault}{\updefault}$A_6$}}}
\put(586,-2041){\makebox(0,0)[lb]{\smash{\SetFigFont{9}{10.8}{\rmdefault}{\mddefault}{\updefault}$w$}}}
\put(5491,-1366){\makebox(0,0)[lb]{\smash{\SetFigFont{9}{10.8}{\rmdefault}{\mddefault}{\updefault}$A_6$}}}
\put(4546,-2041){\makebox(0,0)[lb]{\smash{\SetFigFont{9}{10.8}{\rmdefault}{\mddefault}{\updefault}$A_5$}}}
\put(5041,-2041){\makebox(0,0)[lb]{\smash{\SetFigFont{9}{10.8}{\rmdefault}{\mddefault}{\updefault}$A_4$}}}
\put(6481,-286){\makebox(0,0)[lb]{\smash{\SetFigFont{9}{10.8}{\rmdefault}{\mddefault}{\updefault}$\infty$}}}
\put(2476,-286){\makebox(0,0)[lb]{\smash{\SetFigFont{9}{10.8}{\rmdefault}{\mddefault}{\updefault}$\infty$}}}
\put(6346,-1681){\makebox(0,0)[lb]{\smash{\SetFigFont{9}{10.8}{\rmdefault}{\mddefault}{\updefault} $a_4$}}}
\put(6796,-1816){\makebox(0,0)[lb]{\smash{\SetFigFont{9}{10.8}{\rmdefault}{\mddefault}{\updefault}$a_3$}}}
\put(1936,-1546){\makebox(0,0)[lb]{\smash{\SetFigFont{9}{10.8}{\rmdefault}{\mddefault}{\updefault}$a_2$}}}
\put(5131,-2671){\makebox(0,0)[lb]{\smash{\SetFigFont{9}{10.8}{\rmdefault}{\mddefault}{\updefault}$A_3$}}}
\put(5896,-2671){\makebox(0,0)[lb]{\smash{\SetFigFont{9}{10.8}{\rmdefault}{\mddefault}{\updefault}$A_2$}}}
\put(5761,-1951){\makebox(0,0)[lb]{\smash{\SetFigFont{9}{10.8}{\rmdefault}{\mddefault}{\updefault}$\tilde{w}$}}}
\put(7201,-1411){\makebox(0,0)[lb]{\smash{\SetFigFont{9}{10.8}{\rmdefault}{\mddefault}{\updefault}$a_2$}}}
\end{picture}
$$

\subsubsection*{1.4.3 Combinatorics of the stratification of ${\overline{\cal
      M}}_{0,n+1}(\R)$.}
Note from the proof above that if we relax the condition ``distinct from
the root" in Proposition 1, we get the combinatorics of the stratification of
${\overline{\cal M}}_{0,n+1}(\R)$ instead of its double cover
${\widetilde{\cal M}}_{0,n+1}(\R)$ ($S^{n-2}$ must be replaced by
$\P_{\R}^{n-2}$ in the proof, and cells labeled by $\sigma$ or $\sigma\omega$
must be identified (cf. \S1.2, (1)), where $\omega$ is the permutation $\left(\begin{array}{cccc}
1 & 2 & \ldots & n\\
n & n-1 & \ldots & 1
\end{array}\right)$). In this case, there is no need to use {\it rooted}
trees: strata of ${\overline{\cal M}}_{0,n+1}(\R)$ may now be labeled by
$(n+1)$-unrooted planar trees, $T(\sigma)$, with a cyclic labeling of their
leaves: we mean that the leaves are labeled by $\sigma(0),\ldots,\sigma(n)$, $\sigma\in\Sigma_{n+1}$,
following the trigonometric order of the plane. We may also consider couples
$(T,\sigma)$, where $T$ is an unrooted planar tree, with a given cyclic
labeling of its leaves, from $0$ to $n+1$.\\
Instead of the $\tilde{\nabla}_v$ moves on trees, two distinct moves
$\bar{\nabla}_v^{up}$ and $\bar{\nabla}_v^{down}$ are introduced in
\cite{de}, consisting in reflecting either the ``upper" or the ``lower" subtree
(this distinction is arbitrary)
disconnecting $T$ in $v$. Then the equivalence relation generated by
$T(\sigma)\sim \bar{\nabla}_v^{up}T(\sigma)\sim \bar{\nabla}_v^{down}T(\sigma)$ gives the
combinatorics of ${\overline{\cal M}}_{0,n+1}(\R)$.

\subsection{Translation of the identification rules in terms of nested
  collections and Coxeter groups}

Let $T$ be a planar {\it rooted} $n$-tree, $v$ a vertex of $T$,
$T_v=\{\sigma_i,\sigma_{i+1},\ldots,\sigma_{j-1}\}$ the corresponding subset
of generators of $\Sigma_n$ (cf. \S1.3). Denote by $\omega_{T_v}$ or $\omega _v$ the
 involution $\left(\begin{array}{cccc}
i & i+1 & \ldots & j\\
j & j-1 & \ldots & i
\end{array}\right)$. It is the longest element of $\Sigma_{T_v}$ (the
symmetric group generated by $T_v$) for the word
metric.\\
If ${\cal T}$ is a nested collection (corresponding to an $n$-planar tree
$T$), then for each $T_v$, $T_w$ in ${\cal T}$, define $j_{T_v} T_w$ the connected subset of
$S=\{\sigma_1,\ldots,\sigma_{n-1}\}$ by: $$j_{T_v} T_w=\left\{\begin{array}{cc}
                                          \omega _v T_w \omega _v & \mbox{if } T_w\subset
                                          T_v, \mbox{ or equivalently, }w
                                          \mbox{ is a descendent of } v,\\
                                          T_w & \mbox{if not}.
                                          \end{array} \right.$$ 
\begin{pr}
Under the correspondence $(T,\sigma)\leftrightarrow ({\cal T},\sigma)$,
$(T',\sigma')\leftrightarrow ({\cal T}',\sigma')$ (cf. \S1.3), it holds
$(T',\sigma')={\tilde{\nabla}}_v (T,\sigma)$ if and only if ${\cal T}'=j_{T_v}{\cal T}$ and
$\sigma'=\sigma \omega _{T_v} $, where $j_{T_v}{\cal T}$ is the nested
collection $\{j_{T_v} T_w\in {\cal T}\}$. 
\end{pr}

With the example of figure 2, $\sigma=id$, $\sigma'=\omega_{T_v}=\left(\begin{array}{cccc}
2 & 3 & 4 & 5\\
5 & 4 & 3 & 2
\end{array}\right)$, $T_v=\{\sigma_2,\sigma_3,\sigma_4\}$, $T_w=\{\sigma_2\}$.
The only subset changed under $j_{T_v}$ is  $T_w$: $j_{T_v}T_w=\omega _{T_v}
T_w\omega _{T_v}=\{\sigma_4\}$.\\

A reformulation of the results of \S1.4.2 and \S1.4.3  is\\

\noindent{\bf Theorem A} (Davis-Januszkiewicz-Scott, \cite{d-j-s2}) {\it The two-sheeted covering space $\widetilde{{\cal M}}_{0,n+1}$ is
$\Sigma_n$-equivariantly homeomorphic to the geometric realization $|\Sigma_n {\cal N}|$ of the
following poset $\Sigma_n {\cal N}$:
\begin{itemize}
\item Its elements are the equivalence classes of pairs $({\cal T},\sigma)$, ${\cal
  T}$ a nested collection, $\sigma\in \Sigma_n$,
for the equivalence relation $({\cal T},\sigma)\sim ({\cal T}',\sigma')$ if
and only if there exists  a subset ${\cal T}''\subset{\cal T}$ such that
$\sigma '=\sigma\omega _{{\cal T}''}$ and ${\cal T}'=j_{{\cal T}''}{\cal T}$. Here
$\omega _{{\cal T}''}=\omega _{T_1}\ldots \omega _{T_r}$ if ${\cal T}''=\{T_1,\ldots,T_r\}$,
with $i\leq j$ if $T_i\subset T_j$, and $j_{{\cal T}''}=j_{T_r}\ldots
j_{T_1}$.       
\item The partial order is $[{\cal T},\sigma]\leq [{\cal T}',\sigma']$ if and
  only if there exists some ${\cal T}''$ such that $\sigma '=\sigma\omega _{{\cal
      T}''}$ and $j_{{\cal T}''}{\cal T}\leq {\cal T}'$, i.e. ${\cal T}'\subset j_{{\cal T}''}{\cal T}$.
\end{itemize}
Moreover, the free involution $\tilde{a}:\widetilde{\cal
    M}_{0,n+1}\rightarrow \widetilde{\cal M}_{0,n+1}$ (lifted from the antipodal
    involution of $S^{n-2}$), is combinatorial, and given on the poset
    $\Sigma_n {\cal N}$ by $[{\cal T},\sigma]\mapsto[j_S {\cal
    T},\sigma\omega_S]$, with $S=\{\sigma_1,\ldots,\sigma_{n-1}\}$. Since
    $\overline{\cal M}_{0,n+1}={\widetilde{\cal M}_{0,n+1}}/\tilde{a}$,
    $\overline{\cal M}_{0,n+1}$ inherits a natural cell decomposition with
    poset $\Sigma_n {\cal N}/\tilde{a}$.}

\section{The towers of genus zero real stable curves and action of Tompson and Neretin groups}

\subsection{The towers of genus zero real stable curves ${\overline{\cal
      M}}_{0,3.2^{\infty}}(\R)={\displaystyle\lim_{\stackrel{\longrightarrow}{n}}}\;{\overline{\cal
      M}}_{0,3.2^n}(\R)$ and ${\widetilde{\cal M}}_{0,2^{\infty}}(\R)={\displaystyle\lim_{\stackrel{\longrightarrow}{n}}}\;{\widetilde{\cal
      M}}_{0,2^n +1}(\R)$ }

There is an obvious way to embed ${\overline{\cal M}}_{0,n+1}(\R)$ into
${\overline{\cal M}}_{0,2(n+1)}(\R)$:\\
 If ${\cal C}(x_0,x_1,\ldots,x_n)\in {\overline{\cal
 M}}_{0,n+1}(\R)$ is a stable curve, then graft a new circle at any marked point $x_i$, $i\geq 0$,
with two marked points on it, $y_{2i-1}$ and $y_{2i}$ for $i\not=0$, and
$y_{0}$ and $y_{2n+1}$ for $i=0$. We get a new stable
$2(n+1)$-curve, and it is uniquely defined, since $PGL(2,\R)$ identifies the configurations of any triple of
points on a component.\\
If now ${\cal C}(x_0,x_1,\ldots,x_n)\in {\widetilde{\cal
 M}}_{0,n+1}(\R)$, there are two distinct configurations on the grafted circle with marked
points $y_0, y_{2n+1}$, and double point $x_0$. So we may decide to expand
unambiguously all the points $x_i$ except $x_0$, to get a curve ${\cal
  C}(x_0,y_1,y_2,\ldots, y_{2n-1},y_{2n})$ in ${\widetilde{\cal
 M}}_{0,2n+1}(\R)$.\\
In fact, the map ${\cal C}(x_0,x_1,\ldots,x_n)\mapsto {\cal
  C}(x_0,x_1,\ldots,x_{i-1}, y_i,y_{i+1}, x_{i+1},\ldots,x_n)$ is a section of
a forgetting map ${\overline{\cal M}}_{0,n+2}(\R)\rightarrow {\overline{\cal
    M}}_{0,n+1}(\R)$, and is called ``stabilization" (at least in the complex
case) by Knudsen (cf. \cite{kn}). It is a smooth map.

\begin{pr}[Dyadic expansion principle, or Stabilization]

The em-\\bedding $\overline{\exp}_n: {\overline{\cal M}}_{0,n}(\R)\hookrightarrow
{\overline{\cal M}}_{0,2n}(\R)$ is a morphism of stratified spaces. The
inductive limit  ${\overline{\cal M}}_{0,3.2^{\infty}}(\R)={\displaystyle\lim_{\stackrel{\longrightarrow}{n}}}\;{\overline{\cal M}}_{0,3.2^n}(\R)$ inherits a locally non-finite CW-complex structure.\\
The same is true when the moduli spaces are replaced by their two-sheeted
covering spaces ${\widetilde{\cal M}}_{0,n+1}(\R)$, with embeddings ${\exp_n}: {\widetilde{\cal M}}_{0,n+1}(\R)\hookrightarrow
{\widetilde{\cal M}}_{0,2n+1}(\R)$, and there is a tower
${\widetilde{\cal
    M}}_{0,2^{\infty}}(\R)={\displaystyle\lim_{\stackrel{\longrightarrow}{n}}}\;{\widetilde{\cal M}}_{0,2^n +1}(\R).$\\
We may see both towers as pointed spaces, with based point $*$ represented by
the unique point of ${\overline{\cal M}}_{0,3}(\R)={\widetilde{\cal M}}_{0,2+1}(\R)$. 
\end{pr}        

Proof: We make the proof for the covering spaces ${\widetilde{\cal
    M}}_{0,n+1}(\R)$. If ${\cal M}(T,\sigma)$ is a stratum of ${\overline{\cal M}}_{0,n+1}(\R)$, then
$\exp_n({\cal M}(T,\sigma))$ is the stratum ${\cal M}(\exp_n(T,\sigma))$, where we define $\exp_n(T,\sigma)=(\exp_n(T),\tau)$
to be the planar rooted tree obtained from $(T,\sigma)$ by expanding each of the leaf by two
new edges, which become the input edges of $\exp_n(T,\sigma)$. More precisely, if the input
edges of $T(\sigma)$ are labeled from the left to the right by $\sigma(1),\ldots,
\sigma(n)$, then the input edges of $\exp_n(T)(\tau)$ will be enumerated by $\tau\in
\Sigma_{2n}$ defined by $\tau(2i-1)=2\sigma(i)-1$, $\tau(2i)=2\sigma(i)$,
$i=1,\ldots,n$. By the way, what we have just defined is a group homomorphism
$$\exp_n: \Sigma_n\longrightarrow \Sigma_{2n}:\;\;\sigma\mapsto\tau=\exp_n(\sigma),$$
called the ``expansion morphism", and we have $\exp_n(T,\sigma)=(\exp_n(T),\exp_n(\sigma))$.\\
Clearly, if $(T,\sigma)\sim (T',\sigma')$, then $\exp_n(T,\sigma)\sim \exp_n(T',\sigma')$. Moreover,
though it is not true that  $(T,\sigma)\leq (T',\sigma')\Longrightarrow \exp_n(T,\sigma)\leq \exp_n(T',\sigma')$, this
becomes true by replacing $\exp_n(T,\sigma)$ by an equivalent tree, labeling the same
stratum (using Proposition 1). So $\exp_n$ is cellular. $\;\;\square$\\

%The bracelet model of ${\overline{\cal M}}_{0,n+1}(\R)$ is paticularly
%well-adapted for the description of the embedding $j_n$. Indeed, if
%$b=(h_{\sigma(1)},\ldots,h_{\sigma(n+1)})$ is a labeled bracelet, then
%$j_n(b)$ is the bracelet obtained from $b$ by inserting between $h_{\sigma(i)}$ and
%$h_{\sigma(i+1)}$ the unique horocycle tangent to them. FAIRE UNE FIGURE

\subsection{Action of Thompson and Neretin groups on the towers $\widetilde{\cal M}_{0,2^{\infty}}(\R)$ and $\overline{\cal M}_{0,3.2^{\infty}}(\R)$ }

\subsubsection*{2.2.1 Thompson group $V$ and Neretin group $N$}

{\bf Thompson group.} Recall that Thompson group $V$ (cf. eg. \cite{ca-fl}, \cite{br}) is the group whose elements are represented by symbols
$(\alpha_1,\alpha_0,\sigma)$, where $\alpha_1$ and $\alpha_0$ are binary
planar trees with the same number of leaves, say $n$, and $\sigma\in\Sigma_n$
indicates how the leaves of $\alpha_0$ are mapped to those of $\alpha_1$:
explicitely, if $e_1,\ldots,e_n$ (resp. $f_1,\ldots,f_n$) are the leaves or
input edges of $\alpha_0$ (resp. $\alpha_1$), then $\sigma$ prescribes the
correspondence $e_i\mapsto f_{\sigma(i)}$.\\
 A symbol can be expanded in the
following way: choose an input edge $e_i$, expand its terminal leaf  by two
edges $e_{i,l}$ (the left one), $e_{i,r}$ (the right one), do the same with
$f_{\sigma(i)}$, and define
$\tau\in \Sigma_{n+1}$ describing the correspondence $e_{i,l}\mapsto f_{\sigma(i),l}$,
$e_{i,r}\mapsto f_{\sigma(i),r}$, $e_j\mapsto f_{\sigma(j)}$, $j\not=i$, after
re-labeling the new sets of input edges from $1$ to $n+1$. We get a new symbol 
$(\beta_0,\beta_1,\tau)$, and we consider the equivalence relation generated
by $(\alpha_1,\alpha_0,\sigma)\sim (\beta_0,\beta_1,\tau)$.\\
 The group
structure on the equivalence classes of symbols is defined by

$$(i)\;\;[(\alpha_2,\alpha_1,\tau)][(\alpha_1,\alpha_0,\sigma)]=[(\alpha_2,\alpha_1,\tau\circ\sigma)]$$ 
$$(ii)\;\;[(\alpha_1,\alpha_0,\sigma)]^{-1}=[(\alpha_0,\alpha_1,\sigma^{-1})]$$
$$(iii)\;\;e=[(\alpha,\alpha,id_n)],\; \mbox{for any n-planar tree } \alpha,$$

where we used the fact that a given pair of symbols can be replaced by a pair
of equivalent ones, so that the source tree of the first one coincides with
the target tree of the second one. We have denoted by $e$ the neutral element.\\

Observe that if we expand all the leaves of both trees in a symbol
$(\alpha_1,\alpha_0,\sigma)$, we get a symbol
$(\exp_n(\alpha_1),\exp_n(\alpha_0),\exp_n(\sigma))$, with the notations of the proof
of Proposition 3. So denote the expanded symbol by $\exp_n(\alpha_1,\alpha_0,\sigma)$.\\  

\noindent{\bf Neretin group $N$.} Recall from \cite{ne} or \cite{kap} the
definition of
Neretin ``spheromorphism" group:
Let $\alpha$ be an $n$-planar dyadic rooted tree, and see each leaf of
$\alpha$ as the root of an (infinite) dyadic complete planar tree
$T^{\alpha}_i$, $i=1,\ldots,n$. Thus, $\alpha\cup T^{\alpha}_1\cup\ldots\cup
T^{\alpha}_n$ itself is a dyadic complete planar rooted tree.\\
Now elements of Neretin group $N$ are represented by symbols
$(\alpha_1,\alpha_0,q_{\sigma})$, where again $\alpha_0, \alpha_1$ are
$n$-trees for some $n$, $\sigma\in \Sigma_n$, and $q_{\sigma}$ is a collection
of tree isomorphisms $q_i:T^{\alpha_0}_i\rightarrow T^{\alpha_1}_{\sigma(i)}$,
$i=1,\ldots,n$.\\
The notion of expansion is natural: if $T^{\alpha_0}_i$ is replaced by its two
halves $T^{\alpha_0}_{i,l}$ and $T^{\alpha_0}_{i,r}$, then $q_i$ is replaced
by its restrictions to the halves, $q_{i,l}$ and  $q_{i,r}$, inducing tree
isomorphisms onto the corresponding halves of $T^{\alpha_1}$. But contrary to the Thompson
case, there is no reason that the left (resp. right) half of $T^{\alpha_0}_i$
should be sent onto the left (resp. right) half of $T^{\alpha_1}_i$.\\
The expansion procedure generates an equivalence relation in the set of
symbols, and equivalence classes of symbols composition is induced by the
usual composition of tree isomorphisms. This defines the group structure on
Neretin group $N$.\\
Clearly, $N$ contains $V$ as a subgroup.\\

\noindent{\bf 2.2.2 The cyclic variant.} The use of rooted trees is standard in the definition of Thompson
groups, however we can treat trees occuring in symbols $(\alpha_1,\alpha_0,\sigma)$ as unrooted
$(n+1)$-trees, with $\sigma\in \Sigma_{n+1}$ (instead of $\Sigma_n$). They are
given with a cyclic labeling of their leaves, from $0$ to $n+1$. Now define the notion of
expansion as before. We get a group ${\cal V}$, isomorphic to $V$.\\
An analogous presentation for defining Neretin group is relevant, which uses unrooted
dyadic planar trees: we get a group ${\cal N}$, isomorphic to $N$. In this
way, ${\cal N}$ (and ${\cal V}$) appear as subgroups of the
homeomorphism group of the boundary $\partial {\cal T}_2$ of the regular
dyadic tree ${\cal T}_2$ (without root). In particular, ${\cal N}$ contains as a subgroup the full
automorphism group of the tree, $Aut({\cal T}_2)$, and so, the
p-adic group $PGL(2,\Q_2)$. This illustrates how much bigger ${\cal N}$ is than
${\cal V}$. This corresponds also to the original presentation for ${\cal N}$
(cf. \cite{ne}).\\

\subsubsection*{2.2.3 Thompson and Neretin groups acting on the infinite moduli spaces}

{\bf On the tower ${\widetilde{\cal M}}_{0,2^{\infty}}(\R)$:} Let
$g=[(\alpha_1,\alpha_0,q_{\sigma})]$ be in $ N$, $[{\cal M}(T,\tau)]$ be a stratum of
${\widetilde{\cal M}}_{0,2^{\infty}}(\R)$, represented in some ${\widetilde{\cal
    M}}_{0,2^n+1}(\R)$.\\
At the price of making an expansion of the symbol defining $g$, it can be supposed that the trees $\alpha_i$,
$i=0,1$, are $2^n$-planar trees. Represent the stratum by a symbol $(T_{2^n},T,
\tau)$, where  $T_{2^n}$ is  the dyadic tree:
 $$\stackrel{\underbrace
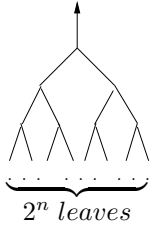
{\begin{picture}(0,0)%
\includegraphics{star.pstex}%
\end{picture}%
\setlength{\unitlength}{2072sp}%
\begingroup\makeatletter\ifx\SetFigFont\undefined%
\gdef\SetFigFont#1#2#3#4#5{%
  \reset@font\fontsize{#1}{#2pt}%
  \fontfamily{#3}\fontseries{#4}\fontshape{#5}%
  \selectfont}%
\fi\endgroup%
\begin{picture}(1689,2127)(979,-1636)
\put(991,-1636){\makebox(0,0)[lb]{\smash{\SetFigFont{6}{7.2}{\rmdefault}{\mddefault}{\updefault}. . .}}}
\put(2296,-1636){\makebox(0,0)[lb]{\smash{\SetFigFont{6}{7.2}{\rmdefault}{\mddefault}{\updefault}. . .}}}
\put(1666,-1636){\makebox(0,0)[lb]{\smash{\SetFigFont{6}{7.2}{\rmdefault}{\mddefault}{\updefault}. . .}}}
\end{picture}}}{2^n\; leaves}$$
It  can be supposed again that $\alpha_0$ and $T_{2^n}$ coincide. Now compose both
symbols in the following way:
$$(\alpha_1,\alpha_0=T_{2^n},q_{\sigma})(T_{2^n},T,\tau)=(\alpha_1,T,\sigma\circ\tau),$$
and making expansions, replace the resulting symbol by an equivalent one, to
get a symbol of the form $(\exp (\alpha_1)=T_{2^m}, \exp(T),\exp(\sigma\circ\tau))$, for some
$m\geq n\in \N$, where $\exp$ is the appropriate expansion map or morphism, and $\exp(\sigma\circ\tau)\in\Sigma_{2^m}$. Finally define
$$g[{\cal M}(T,\tau)]:=[{\cal M}(\exp (T),\exp(\sigma\circ\tau))].$$

\noindent{\bf On the cyclic tower ${\overline{\cal M}}_{0,3.2^{\infty}}(\R)$:}
We now use the
presentation of ${\cal N}$ (and ${\cal V}$) given in the variant above, to make symbols
defining group elements act on symbols labeling strata. As can be guessed, the
tree $T_{2^n}$ must be replaced here by the unrooted tree $T_{3.2^n}$. Trees
involved in symbols defining elements of ${\cal N}$ or ${\cal V}$ must be represented with a marked vertex, corresponding to the central
vertex of the tree $T_{3.2^n}$. If now $\overline{\cal M}(T(\tau))$ is a
stratum of $\overline{\cal M}_{0,3.2^n}(\R)$,
$g=[(\alpha_1,\alpha_0=T_{3.2^n},q_{\sigma})]\in {\cal N}$, with $\sigma\in
\Sigma_{3.2^n}$, then compose the symbols:
$(\alpha_1,\alpha_0=T_{3.2^n},q_{\sigma}).(T_{3.2^n},T(\tau))=(\alpha_1,T(\sigma\circ\tau))$,
and expand $\alpha_1$ to get a tree of the form
$T_{3.2^m}=\overline{\exp}(\alpha_1)$ (the marked
vertex is useful to perform this procedure), and correspondingly expand
$T(\sigma\circ\tau)$. Then define $g.[\overline{\cal M}(T(\tau))]=[\overline{\cal M}(\overline{\exp}(T(\sigma\circ\tau))]$.    
 
\begin{theo}
Denoting by $Strat\left({\widetilde{\cal M}}_{0,2^{\infty}}(\R)\right)$ the set of
strata of the tower ${\widetilde{\cal M}}_{0,2^{\infty}}\\(\R)$, the map $$N\times
Strat\left({\widetilde{\cal M}}_{0,2^{\infty}}(\R)\right) \longrightarrow Strat\left({\widetilde{\cal M}}_{0,2^{\infty}}(\R)\right)$$
defined above is well-defined, and induces a cellular left
action of Neretin group $N$ on the tower ${\widetilde{\cal
    M}}_{0,2^{\infty}}(\R)$:
$$\tilde{\gamma}: N\longrightarrow Homeo_{cell}({\widetilde{\cal M}}_{0,2^{\infty}}(\R)),
$$ which is faithful.\\
In particular, $\tilde{\gamma}$ restricts to a cellular action of Thompson group
$V$.\\
There is an analogous statement with $\widetilde{\cal M}_{0,2^{\infty}}(\R)$
replaced by $\overline{\cal M}_{0,3.2^{\infty}}(\R)$, and $\tilde{\gamma}$ by
$\overline{\gamma}: {\cal N}\longrightarrow Homeo_{cell}({\overline{\cal M}}_{0,3.2^{\infty}}(\R))$.         
\end{theo}

Proof: We prove the statement for the action on $\widetilde{\cal
  M}_{0,2^{\infty}}(\R)$:\\
(1) We first need to check  $$(T,\tau)\sim
(T',\tau ')\Longrightarrow
(\exp(T),\exp(\sigma\circ\tau))\sim
(\exp(T'),\exp(\sigma\circ\tau ')):$$
It is clear we can restrict to the case where ${\cal T}'=j_{T_v}{\cal T}$ and
$\tau '= \tau \omega _v$, for some $v\in Vert(T)$ (we use the correspondence
between planar trees and nested families, cf. \S1.3). For simplicity suppose $\omega _v$ is of the form  $$\omega _v= 
 \left(\begin{array}{cccc}
1 & 2 & \ldots & k\\
k & k-1 & \ldots & 1
\end{array}\right)=\omega _{(1,\ldots,k)}$$
Using the
well-defined homorphism $\exp:\Sigma_{2^n}\rightarrow \Sigma_{2^m}$ it follows
that $\exp(\sigma\circ\tau ')=\exp(\sigma\circ\tau)\exp(\omega _v)$. Suppose
for simplicity again that $m=n+1$; then with the labeling of $\exp(T)$, it holds
$$\exp(\omega _v)= \left(\begin{array}{cccccccc}
1 & 2 & \ldots & 2l-1 & 2l & \ldots & 2k-1 & 2k\\
2k-1 & 2k & \ldots & 2k-2l+1 & 2k-2l+2 & \ldots & 1 & 2
\end{array}\right).$$ 
Denoting by $\tilde{v}$ the vertex $v$ seen in the expanded tree $\exp(T)$, $\omega _v$ differs from $$\omega _{\tilde v}= \left(\begin{array}{cccccccc}
1 & 2 & \ldots & 2l-1 & 2l & \ldots & 2k-1 & 2k\\
2k & 2k-1 & \ldots & 2k-2l+2 & 2k-2l+1 & \ldots & 2 & 1
\end{array}\right)=\omega _{(1,\ldots,2k)}$$ by the product of the $k$ transpositions
$\omega _{\tilde{v}_1}\ldots\omega _{\tilde{v}_k}$, where $\omega
_{\tilde{v}_i}$ is simply $\omega_i=(2i-1,2i)$, and $\tilde{v}_i$ is the $i^{th}$ leaf
of $T$ but seen in $\exp(T)$: each $\tilde{v}_i$ has itself two descendent
leaves. It follows that $$\left\{\begin{array}{c}

\exp({\cal T}')=j_{T_{\tilde{v}_1}}\ldots j_{T_{\tilde{v}_k}} j_{T_{\tilde{v}}}\exp({\cal T}),\\
\exp(\sigma\circ\tau ')= \exp(\sigma \circ \tau)\omega _{\tilde{v}}\omega _{\tilde{v}_k}\ldots\omega _{\tilde{v}_1},
\end{array}\right.$$
where in fact the operations $j_{T_{\tilde{v}_i}}$ have no effect. The
expected equivalence is proved.\\

(2) We must check also that
$$\exp_n (\alpha_1,\alpha_0=T_{2^n},q_{\sigma})\exp_n (T_{2^n},T,\tau) \sim \exp_n
(\alpha_1,T,\sigma\circ\tau).$$
Now $\exp_n
(\alpha_1,T_{2^n},q_{\sigma})=(\exp_n(\alpha_1),T_{2^{n+1}},\tilde{q}_{\tilde{\sigma}}),$ where $\tilde{\sigma}$ may differ from
$\exp_n(\sigma)$ (because $g$ sits in $N$, not necessarily in $V$) by a product of transpositions $\omega_i=(2i-1,2i)$:
$\tilde{\sigma}=\exp(\sigma)\omega_{i_1}\ldots\omega_{i_r}$. So,
$$(*)\;\; \exp_n(\alpha_1,T_{2^n},q_{\sigma})\exp_n (T_{2^n},T,\tau)=(\exp_n(\alpha_1),
\exp_n(T),\exp_n(\sigma)\omega_{i_1}\ldots$$
$$\ldots\omega_{i_r}\exp_n(\tau)).$$ 
But $\exp(\tau)\omega_i=\omega_{\tau(i)}\exp(\tau)$, so that
$$(*)=(\exp_n(\alpha_1),\exp_n(T),\exp_n(\sigma)\exp_n(\tau)\omega_{\tau^{-1}(i_1)}\ldots\omega_{\tau^{-1}(i_r)})$$
$$\sim (\exp_n(\alpha_1),\exp_n(T),\exp_n(\sigma\circ\tau)).$$

(3) Faithfulness is obvious, or may be proven by invoking the simplicity of the
group $N$ (\cite{kap}).\\

(4) Translated in terms of $\tilde{\nabla}_v$ moves, the proof above can be now easily adapted to the case of the action on the tower
    $\overline{\cal M}_{0,3.2^{\infty}}(\R)$ (but using $\bar{\nabla}_v^{up}$
    and $\bar{\nabla}_v^{down}$ moves). $\;\;\square$\\

{\bf Two complementary remarks:}\\

{\bf 1. Pointwise description of the action.} Let $(a_i)_{i=1,\ldots,n}$ be the
affine coordinates of a point $p\in {\cal M}(T,\tau)$ (cf. proof of
Proposition 1,
where we give the coordinates for a codimension 1 stratum). If $v$ is an
internal vertex of $T$, and $i_1<\ldots <i_k$ is the collection of the
{\it one-valent} descendents (``sons") of $v$ (recall the labeling of $T$ is canonical, from
$1$ (on the left) to $n$ (on the right)), this means that  
$$a_{\tau(i_1)}<\ldots <a_{\tau(i_k)}.$$ 
Then $g(p)$, with $g=[(\alpha_1,\alpha_0,q_{\sigma})]$ in $N$, corresponds to
an adequate expansion in  ${\cal M}(\exp(T),\exp(\sigma\circ\tau))$, of the
point of ${\cal M}(T,\sigma\circ\tau)$ with coordinates $(b_i)_{i=1,\ldots,n}$ such
that $b_{\sigma(i)}=a_i$.\\ 

{\bf 2. Action on the complex tower.} As was pointed out to the author by
Yu. Neretin, we could define also the pointwise action of $N$ on the complex tower
$\overline{\cal M}_{0,2^{\infty}}(\C)$ by the same formula as above, replacing ${\cal
  M}(T,\tau)$ by the complex stratum ${\cal M}(T)$ (the strata are indeed
labeled by bare $n$-rooted trees, isotopically realized in $\R^3$) and the real coordinates by complex ones. This
is an action of stratified spaces, and Theorem 1 shows that the restriction of
the action to the real part respects the cell decomposition in
associahedra of the latter. We mention also that in the cyclic version of the action, the complex
action of ${\cal N}$ extends the action of $Aut({\cal T}_2)$ that was observed
by Kapranov in the note \cite{ka3}.

\subsection{Thompson group ${\cal T}$ acting on the cyclic tower $\overline{\cal
    M}_{0,3.2^{\infty}}(\R)$}

Thompson group ${\cal T}$ is the subgroup of elements of ${\cal V}$ represented by symbols  
$(\alpha_1,\alpha_0,\sigma)$ (in the second description) where $\sigma$ lies
in the cyclic subgroup $\Z/(n+1)\Z$ of $\Sigma_{n+1}$ (when $\alpha_1$ and
$\alpha_0$ are unrooted $(n+1)$-trees).\\
On the other hand, since the expansion morphism
$\overline{\exp}_n:\Sigma_n\rightarrow \Sigma_{2n}$ maps $\Z/n\Z$ into $\Z/2n\Z$,
there is a distinguished subcomplex in the tower $\overline{\cal
  M}_{0,3.2^{\infty}}(\R)$, which is the union of strata ${\cal M}(T(\tau))$
with $\tau$ in the cyclic group $\Z/n\Z$ (if $T(\tau)$ has $n$ leaves). This
subcomplex is simply an infinite associahedron $K_{\infty}$, an inductive
limit of associahedra $K_{3.2^n -1}$. The cyclic tower $\overline{\cal M}_{0,3.2^{\infty}}(\R)$ is tiled by an
infinite number of copies of the associahedron $K_{\infty}$. 

\begin{theo}
The stabilizer of the distinguished associahedron $K_{\infty}$ under the action of
${\cal N}$ is the semi-direct product $Out({\cal T})\triangleright {\cal T}$, where $Out({\cal T})$, the exterior
automorphism group of ${\cal T}$, is known to be the order two cyclic group. 
\end{theo}

 Remark: $Out({\cal T})\triangleright {\cal T}$ as a symmmetry group of $K_{\infty}$ is
 the infinite-dimensional analogue of $D_n=\Z/2\Z\triangleright\Z/n\Z$,
 the dihedral group, as the symmmetry group of the associahedron $K_{n-1}$.\\
 
Proof: A cell of $K_{\infty}$ is represented in some $\overline{\cal
  M}_{0,3.2^n}(\R)$ by some ${\cal M}(T(\tau))$, with $\tau\in \Z/3.2^n\Z$. If
  $g\in {\cal T}$ is represented by $(\alpha_1,\alpha_0=T_{3.2^n},\sigma)$, with $\sigma\in
  \Z/3.2^n\Z$, then $g.[{\cal M}(T(\tau))]=[{\cal M}(\exp(T(\sigma\circ\tau)))]$,
  where $\sigma\circ\tau\in \Z/3.2^n\Z$, and so $\exp(T(\sigma\circ\tau))$ is a
  tree $T'(\sigma')$, with $\sigma'$ in some $\Z/3.2^m\Z$. This proves that
  $g.[{\cal M}(T(\tau))]$ is a cell of $K_{\infty}$.\\
Now, if $g=[(\alpha_1,\alpha_0=T_{3.2^n},q_{\sigma_0})]$ lies in ${\cal N}$ and
  stabilizes $K_{\infty}$, then at the price of composing $g$ by an element of
  $T$ it may be supposed that $\alpha_1=\alpha_0= T_{3.2^n}$, and that
  $\sigma_0(0)=0$. Now the  the maximal cell of $\overline{\cal
  M}_{0,3.2^n}(\R)\cap  K_{\infty}$, labeled by the stared tree $T(id_{3.2^n})$, is
  mapped by $g$ to the cell ${\cal M}(\exp(T(\sigma_0)))$. Using
  $\bar{\nabla}$ moves, we see that this cell still lies
  in $K_{\infty}$ if and only if $\sigma_0=id_{3.2^n}$ or
  $\sigma_0=\omega_{(1,\ldots,3.2^n)}$ (using Theorem A). Repeating the argument in $\overline{\cal
  M}_{0,3.2^{n+1}}(\R)$ with $\sigma_1$ instead of $\sigma_0$, it comes
  $\sigma_1=id_{3.2^{n+1}}$ or $\sigma_1=\omega_{(1,\ldots,3.2^{n+1})}$. But
  $\sigma_1$ is the expansion of $\sigma_0$, possibly twisted by elementary
  transpositions, so
  $\sigma_0=id_{3.2^n}\Longrightarrow\sigma_1=id_{3.2^{n+1}}$ and
  $\sigma_0=\omega_{(1,\ldots,3.2^n)}\Longrightarrow \sigma_1=
  \omega_{(1,\ldots,3.2^{n+1})}$. Then we see that $g=id$ or the involution $Inv$
  defined by $\sigma_0=\omega_{(1,\ldots,3.2^n)}$,
  $\sigma_1=\omega_{(1,\ldots,3.2^{n+1})}$, etc. Under the isomorphism $T\cong
  PPSL(2,\Z)$ (cf. \cite{im}), $Inv$ corresponds to the non-preserving orientation inversion
  $x\in \R P^1\mapsto \frac{1}{x}\in \R P^1$. And it is known to generate the
  exterior automorphism group of ${\cal T}$, $Out({\cal T})\cong \Z/2\Z$ (cf. \cite{bri}). $\;\;\square$\\               
 
%These actions are far from being free, as shown in the next lemma:

%\begin{lem}
%The stabilizer of the base point of ${\widetilde{\cal M}}_{0,2^{\infty}}(\R)$ under the
%action of Neretin group is
%isomorphic to the automorphism group of the dyadic rooted infinite tree. As
%for the stabilizer of the base point of ${\overline{\cal
%    M}}_{0,3.2^{\infty}}(\R)$, it is isomorphic to the stabilizer of a vertex
%in the automorphism group $Aut({\cal T}_2)$ of the (unrooted) dyadic complete tree.
%\end{lem}

%Proof: We treat the first case: let $*$ be the 0-cell corresponding to the point $\widetilde{\cal M}_{0,3}(\R)$, represented by
%the symbol $(T_{2^n},T_{2^n}, id_{2^n})$, and $g\in N$ represented by
%$(\alpha_1,T_{2^n},q_{\sigma})$, where $\alpha_1$ is a $2^n$-tree and $\sigma\in
%\Sigma_{2^n}$. So $g.*=*\iff (\alpha_1,T_{2^n},\sigma)\sim
%(T_{2^n},T_{2^n}, id_{2^n})$. Since $\alpha_1$ and $T_{2^n}$ are both
%$2^n$-labeled trees, the only possibility is $\alpha_1=T_{2^n}$, and $\sigma\in
%\Sigma_{2^n}$ must belong to the automorphism group $Aut(T_{2^n})$. The
%conclusion follows easily. $\;\;\square$

\subsection{Extensions of Thompson and Neretin groups by infinite pure
  quasi-braid groups $PJ_{2^{\infty}}$ and $Q_{3.2^{\infty}}$}

It is shown in \cite{d-j-s1} that ${\overline{\cal M}}_{0,n+1}(\R)$ (or ${\widetilde{\cal M}}_{0,n+1}(\R)$) is an
aspherical space: its universal covering ${\widehat{\cal M}}_{0,n+1}(\R)$
is contractible, so that its cellular homology coincides with the group homology of its
fundamental group. The universal covering of ${\widetilde{\cal
    M}}_{0,2^{\infty}}(\R)$ is the inductive limit of the coverings
${\widehat{\cal M}}_{0,2^n+1}(\R)$, and will be denoted ${\widehat{\cal M}}_{0,2^{\infty}}(\R)$. Let us denote by $PJ_n$ the fundamental group of
${\widetilde{\cal M}}_{0,n+1}(\R)$: $PJ_n=\pi_1({\widetilde{\cal
    M}}_{0,n+1}(\R))$. Since ${\widetilde{\cal M}}_{0,n+1}(\R)$ is a
double-cover of ${\overline{\cal M}}_{0,n+1}(\R)$, there is an extension
$1\rightarrow PJ_n\rightarrow Q_{n+1}\rightarrow \Z/2\Z\rightarrow 0$, where
we set $Q_{n+1}:=\pi_1({\widetilde{\cal M}}_{0,n+1}(\R))$.\\

Let us now consider the group $PJ_{2^{\infty}}=\pi_1({\widetilde{\cal
    M}}_{0,2^{\infty}}(\R))$, so that   
$$PJ_{2^{\infty}}={\displaystyle\lim_{\stackrel{\longrightarrow}{n}}
  PJ_{2^n}},$$
as well as the group $Q_{3.2^{\infty}}=\pi_1({\overline{\cal
    M}}_{0,3.2^{\infty}}(\R))={\displaystyle\lim_{\stackrel{\longrightarrow}{n}}  Q_{3.2^n}}.$ Now each transformation $\tilde{\gamma}(g):{\widetilde{\cal
    M}}_{0,2^{\infty}}(\R)\longrightarrow {\widetilde{\cal M}}_{0,2^{\infty}}(\R)$,
with $g\in N$, can be lifted -- non-uniquely -- to the universal covering
${\widehat{\cal M}}_{0,2^{\infty}}(\R)$, and similarly for the transformation $\bar{\gamma}(g)$.

\begin{de}
Let $G$ denote Thompson group $V$ or Neretin group $N$. The group of lifted transformations $\tilde{\gamma}(g)$
(resp. $\bar{\gamma}(g)$) G is denoted $\tilde{\cal A}_G$ (resp. $\bar{\cal A}_G$), and is a subgroup
of the cellular homeomorphism group of ${\widehat{\cal M}}_{0,2^{\infty}}(\R)$
(resp. ${\widehat{\cal M}}_{0,3.2^{\infty}}(\R)$). 
\end{de}

By construction we get epimorphisms $\tilde{\cal A}_G\longrightarrow G$ and
$\bar{\cal A}_{\cal G}\longrightarrow {\cal G}$. The kernel of the first epimorphism is the
automorphism of the universal covering map ${\widehat{\cal
    M}}_{0,2^{\infty}}(\R)\longrightarrow {\widetilde{\cal
    M}}_{0,2^{\infty}}(\R)$: it coincides with the fundamental group
$PJ_{2^{\infty}}$. Similarly, the kernel of the second one
$Q_{3.2^{\infty}}$. We are now ready for introducing the announced extensions
of Thompson and Neretin groups:

\begin{de}[Quasi-braid extensions of Thompson and Neretin groups]

Let $G$ (resp. ${\cal G}$)  denote Thompson group $V$ (resp. ${\cal V}$) or
Neretin group $N$ (resp. ${\cal N}$). The construction above
provides extensions $1\rightarrow PJ_{2^{\infty}}\longrightarrow \tilde{\cal A}_G \longrightarrow
G\rightarrow 1$ and $1\rightarrow Q_{3.2^{\infty}}\longrightarrow \bar{\cal
  A}_{\cal G} \longrightarrow {\cal G}\rightarrow 1$, and the natural embedding
$G\hookrightarrow {\cal G}$ induces a morphism between both extensions
$$1\rightarrow PJ_{2^{\infty}}\longrightarrow \tilde{\cal A}_G \longrightarrow
G\rightarrow 1,$$
$$ 1\rightarrow Q_{3.2^{\infty}}\longrightarrow \bar{\cal A}_{\cal G} \longrightarrow {\cal G}\rightarrow
1,$$
where all the vertical arrows are injective morphisms.
\end{de}

The embedding $G\hookrightarrow {\cal G}$ is induced by the embeddings of the
rooted planar trees $T_{2^n}$ as subtrees of the unrooted planar tree
$T_{3.2^n}$. The latter induce also the cellular embeddings $\overline{\cal M}_{0,2^n +1}\hookrightarrow \overline{\cal M}_{0,3.2^n}$, and morphisms
$Q_{2^n +1}\hookrightarrow Q_{3.2^n}$ (they are injective, cf. Theorem
4). Composing $PJ_{2^n}\hookrightarrow Q_{2^n +1}$ with $Q_{2^n
  +1}\hookrightarrow Q_{3.2^n}$ gives $PJ_{2^n}\hookrightarrow Q_{3.2^n}$, and
so $PJ_{2^{\infty}}\hookrightarrow Q_{3.2^{\infty}}$.

\section{An analogue of the Euler class for Neretin Spheromorphism group}

Our first task is to improve our understanding of the infinite quasi-braid groups.

\subsection{Quasi-braid groups $J_n$}
{\bf 3.1.1} Let $J_n$ be the group defined in \cite{d-j-s2} by generators and relations, with generators
$\alpha_T$ for each subset of $S=\{\sigma_1,\ldots,\sigma_{n-1}\}$ such that
the corresponding graph $G_T$ is connected, with the following relations:
\begin{itemize}
\item $\alpha_T ^2=1$ for each $T$
\item $\alpha_T \alpha_{T'}=\alpha_{j_T T'}\alpha_T$ if $T'\subset T$
\item $\alpha_T \alpha_{T'}=\alpha_{T'} \alpha_T$ if $G_{T\cup T'}$ is not
  connected.
\end{itemize}

These relations are those verified by the involutions $\omega_T$. So there is
a well-defined homomorphism $\phi: J_n\rightarrow \Sigma_n$, $\alpha_T\mapsto \omega_T$,
surjective since for $T_i=\{\sigma_i=\\(i,i+1)\}$, $\phi(\alpha_{T_i})=\sigma_i$.\\

\noindent{\bf Theorem B} (Davis-Januszkiewicz-Scott, \cite{d-j-s2}) {\it The universal cover $\widehat{{\cal M}}_{0,n+1}$ of the two-sheeted cover $\widetilde{{\cal M}}_{0,n+1}$ is
$J_n$-equivariantly homeomorphic to the geometric realization $|J_n {\cal N}|$
of the poset $J_n {\cal N}$:
\begin{itemize}
\item Its elements are the equivalence classes of pairs $({\cal T},\alpha)$, ${\cal
  T}$ a nested collection, $\alpha\in J_n$,
with the equivalence relation $({\cal T},\alpha)\sim ({\cal T}',\alpha')$ if
and only if there exists  a subset ${\cal T}''\subset{\cal T}$ such that
$\alpha '=\alpha\alpha _{{\cal T}''}$ and ${\cal T}'=j_{{\cal T}''}{\cal T}$. Here
$\alpha _{{\cal T}''}=\alpha _{T_1}\ldots \alpha _{T_r}$ if ${\cal T}''=\{T_1,\ldots,T_r\}$,
with $i\leq j$ if $T_i\subset T_j$, and $j_{{\cal T}''}=j_{T_r}\ldots
j_{T_1}$.       
\item The partial order is $[{\cal T},\alpha]\leq [{\cal T}',\alpha']$ if and
  only if there exists some ${\cal T}''$ such that $\alpha '=\alpha\alpha _{{\cal
      T}''}$ and $j_{{\cal T}''}{\cal T}\leq {\cal T}'$, i.e. ${\cal
    T}'\subset j_{{\cal T}''}{\cal T}$.
\end{itemize}
Moreover, there is a natural $J_n$-left-equivariant map $J_n {\cal
  N}\rightarrow \Sigma_n {\cal N}$  given by $[{\cal
  T},\alpha]\rightarrow [{\cal T},\phi(\alpha)]$, the $J_n$-action on
$\Sigma_n {\cal N}$ being defined by $\alpha: [{\cal T},\sigma]\mapsto [{\cal
  T},\phi(\alpha)\sigma]$.\\
The kernel $PJ_n:=Ker\phi$ is the fundamental group of $\widetilde{{\cal
    M}}_{0,n+1}$, and there is a short exact sequence
$$1\rightarrow PJ_n \longrightarrow J_n \longrightarrow \Sigma_n \rightarrow
1.$$}

As for the fundamental group of $\overline{{\cal M}}_{0,n+1}$, it is generated
by $Ker\; \phi=PJ_n$ and a lift $\hat{a}$ of the free involution map
$\tilde{a}:[{\cal T},\sigma]\mapsto [j_S{\cal T},\sigma\omega_S]$ (cf. Theorem
A), that we now describe:\\ 

-- Let $\hat{\alpha}:=\alpha_{T_1}\ldots\alpha_{T_r}$ be any lift in $J_n$ of
 $\omega_S$. We choose the unique lift $\hat{a}$ such that
 $\hat{a}[\emptyset, 1]=[\emptyset, \hat{\alpha}]$.\\

-- Let $[\emptyset,\beta]$ be a maximal cell of $|J_n {\cal N}|$, with
 $\beta=\alpha_{U_1}\ldots\alpha_{U_k}$. Choose an edge path $[\emptyset,
 1]\rightarrow [\emptyset,\alpha_{U_1} ]\rightarrow
 [\emptyset,\alpha_{U_1}\alpha_{U_2}]\rightarrow\ldots\rightarrow[\emptyset,\alpha_{U_1}\ldots\alpha_{U_k}=\beta]$
 between $[\emptyset, 1]$  and $[\emptyset,\beta]$ (it is an edge path in the
 dual complex of $|J_n {\cal N}|$). Then
$$\hat{a}[\emptyset,\beta]=[\emptyset,\hat{\alpha}\alpha_{j_S
  {U_1}}\ldots\alpha_{j_S U_k}].$$
It is easy to check that the correspondence $j_S: \alpha_T\mapsto \alpha_{j_S T}$
induces a well-defined automorphism of the group $J_n$ (so that $j_S \beta=\alpha_{j_S  {U_1}}\ldots\alpha_{j_S U_k}$  does not depend on the way of writing $\beta$). Now $\hat{a}^2$ is the
map such hat 
$$\hat{a}^2[\emptyset,\beta]=[\emptyset,\hat{\alpha} j_S(\hat{\alpha})\alpha_{U_1}\ldots\alpha_{U_k}].$$
In other words, $\hat{a}^2$ is the pure quasi-braid $\hat{\alpha}j_S(\hat{\alpha})=\alpha_{T_1}\ldots\alpha_{T_r}\alpha_{j_S {T_1}}\ldots\alpha_{j_S
  {T_r}} $, acting on the left on the cells. This was of course predictable
from the short exact sequence $1\rightarrow PJ_n\rightarrow Q_{n+1}\rightarrow
\Z/2\Z\rightarrow 1$.\\              

\noindent{\bf 3.1.2 Describing the morphism $PJ_n\rightarrow PJ_{2n}$ and
  defining $J_n\rightarrow J_{2n}$}. The formalism contained in Theorem B enables us to describe the morphism $PJ_n\rightarrow
PJ_{2n}$ induced by the embedding $\exp_n:\widetilde{{\cal
    M}}_{0,n+1}\hookrightarrow \widetilde{{\cal M}}_{0,2n+1}$:
Each $\alpha=\alpha_{T_1}\ldots\alpha_{T_r}$ in $PJ_n$ projects onto
$\omega_{T_1}\ldots\omega_{T_r}=1$ in $\Sigma_n$. We interpret $\alpha$ as the
homotopy class of the
edge loop $\gamma=(id_n\rightarrow
\omega_{T_1}\rightarrow\omega_{T_1}\omega_{T_2}\rightarrow\ldots\rightarrow\omega_{T_1}\omega_{T_2}\ldots\omega_{T_r}=id_n)$
in the dual cell complex of $|\Sigma_n {\cal N}|\cong \widetilde{{\cal M}}_{0,n+1}$.\\
By definition of the embedding $\exp_n:\widetilde{{\cal
    M}}_{0,n+1}\hookrightarrow \widetilde{{\cal M}}_{0,2n+1}$, the loop
$\gamma$ is mapped to the loop
 $\exp(\gamma)=(id_{2n}\rightarrow\exp(\omega_{T_1})\rightarrow\exp(\omega_{T_1})\exp(\omega_{T_2})\rightarrow\ldots\rightarrow\exp(\omega_{T_1})\exp(\omega_{T_2})\ldots\exp(\omega_{T_r})=id_{2n}).$ We need however to
 precise what a path of the form $id_{2n}\rightarrow\exp(\omega_{T})$ is:\\
Suppose for simplicity that $\omega_T$ is of the form  $$\omega _T= 
 \left(\begin{array}{cccc}
1 & 2 & \ldots & k\\
k & k-1 & \ldots & 1
\end{array}\right),$$ so that  
$\exp(\omega_{T})=\left(\begin{array}{ccccc}
1 & 2 & \ldots & 2k-1 & 2k \\
2k-1 & 2k & \ldots & 1 & 2
\end{array}\right)$ is the product $\omega_{\exp(T)}\omega_{(1,2)}\\\ldots \omega_{(2k-1,2k)}$,
with $\omega_{(2i-1,2i)}$ the transposition $\sigma_{2i-1}$. The path
 $id_n\rightarrow \omega_{T}$ once embedded in $|\Sigma_{2n}{\cal N}^0|$, and
 after a suitable translation to make its extremities coincide with the
 barycenters of the cells $id_{2n}$ and $\exp(\omega_T)$ (which are adjacent
 because $\exp$ is cellular, meeting along a codimension $n+1$ stratum), becomes the
 straight line joining the barycenters. We claim this line is homotopic to the
 edge path
 $$id_{2n}\rightarrow\omega_{\exp(T)}\rightarrow\omega_{\exp(T)}\omega_{(1,2)}\rightarrow\ldots\rightarrow\omega_{\exp(T)}\omega_{(1,2)}\ldots\omega_{(2k-1,2k)}.$$
Indeed, the path above passes through cells which all share a same codimension
$n+1$ stratum, and the line $id_{2n}\rightarrow \exp(\omega_T)$ crosses the
same stratum.\\ 

  Now $\alpha_T$ may be lifted in $J_{2n}$ to
$$\exp(\alpha_T):= \alpha_{\exp(T)}\alpha_{(1,2)}\ldots \alpha_{(2k-1,2k)},$$ where
$\alpha_{(2i-1,2i)}:=\alpha_{T_i}$, with $T_i=\{\sigma_{2i-1}\}$. Finally define
$\exp(\alpha)$ as the product
$\exp(\alpha_{T_1})\ldots\exp(\alpha_{T_r})$. We now claim:  

\begin{pr}
\begin{enumerate}
\item The map $J_{n}\rightarrow J_{2n}$ : $\alpha\mapsto \exp(\alpha)$, is a
well-defined group homomorphism. More generally, each dyadic expansion map
(not morphism)
$\Sigma_{n}\rightarrow \Sigma_{n+*}$ has a canonical lift $J_{n}\rightarrow J_{n+*}$.
\item Its restriction to $PJ_n$ is the morphism $(\exp_n)_*: PJ_{n}\rightarrow PJ_{2n}$ induced at the fundamental group level by
the embedding $\exp_n:\widetilde{{\cal M}}_{0,n+1}\hookrightarrow \widetilde{{\cal M}}_{0,2n+1}$.             
\item The morphisms $(\exp_n)_*$ are injective for all $n\geq 2$. In particular,
  the group $PJ_{2^{\infty}}={\displaystyle\lim_{\stackrel{\longrightarrow}{n}}
    PJ_{2^n}}$ is the limit of an inductive system of embeddings. The same
  injectivity property holds for the homology maps
  $(\exp_n)_*:H_*(PJ_n,\Z)\rightarrow H_*(PJ_{2n},\Z)$.
\item There are commutative diagrams\\
\setlength{\unitlength}{0.9cm}

\begin{picture}(10,2) 
\multiput(4,2)(1.5,0){4}{\vector(1,0){0.5}}   %les vecteurs horizontaux en haut
\put(3.5,1.9){1} \put(4.6,1.9){$PJ_n$} \put(6.3,1.9){$J_n$}
\put(7.8,1.9){$\Sigma_n$}
\put(9.2,1.9){1}

\put(6.5,1.4){\vector(0,-1){.6}}  %le vecteur vertical 
\multiput(4,0.4)(1.5,0){4}{\vector(1,0){0.5}}   %les vecteurs horizontaux en bas
\put(3.5,0.3){1} \put(4.6,0.3){$PJ_{2n}$} 
\put(6.3,0.3){$J_{2n}$}
\put(7.8,0.3){$\Sigma_{2n}$}
\put(9.2,0.3){1}

%\put(8.5,3.4){\line(0,-1){.7}}\put(8.6,3.4){\line(0,-1){.7}}
\put(4.9,1.4){\vector(0,-1){.6}}

%\put(10.2,4){\vector(1,0){1}} \put(10.2,2.1){\vector(1,0){1}}
\put(8,1.4){\vector(0,-1){.6}}
\end{picture}

where all the vertical maps are injective morphisms, producing the short exact
sequence
$$1\rightarrow PJ_{2^{\infty}}\longrightarrow J_{2^{\infty}} \longrightarrow
\Sigma_{2^{\infty}}\rightarrow 1,$$ where
$J_{2^{\infty}}={\displaystyle\lim_{\stackrel{\longrightarrow}{n}} J_{2^n}}$
and $\Sigma_{2^{\infty}}={\displaystyle\lim_{\stackrel{\longrightarrow}{n}}
  \Sigma_{2^n}}$. It may be obtained by restricting the extension $1\rightarrow PJ_{2^{\infty}}\longrightarrow \tilde{\cal A}_{V} \longrightarrow
V\rightarrow 1$ to the subgroup $\Sigma_{2^{\infty}}\subset V$. 
\end{enumerate}
\end{pr}

Proof: 1) We must check that $\exp$ preserves the relations of the group
$J_n$:
For simplicity, suppose $G_T=(1,\ldots,j)$ and compute $\exp(\alpha_T)^2$:   
 $\exp(\alpha_T)^2=\alpha_{\exp(T)}\alpha_{(1,2)}\ldots
 \alpha_{(2j-1,2j)}\alpha_{\exp(T)}\alpha_{(1,2)}\ldots \alpha_{(2j-1,2j)}$. Now
 oberve that in $J_{2n}$, for all $i\leq j$,
 $\alpha_{\exp(T)}\alpha_{(2i-1,2i)}=\alpha_{2(j-i+1)-1,2(j-i+1)}\alpha_{\exp(T)}$. This fact joint to the commutation property of $\alpha_{(1,2)},\ldots ,\alpha_{(2j-1,2j)}$ among themselves allows to write $\exp(\alpha_T)=\alpha_{(1,2)}\ldots \alpha_{(2j-1,2j)}\alpha_{\exp(T)}$, and it comes easily $\exp(\alpha_T)^2=1$.\\
Then let $T'\subset T$ be such that $G_{T'}=(1,\ldots,i)$, $i\leq j$, and
check the relation $\alpha_T \alpha_{T'}\alpha_T=\alpha_{j_T T'}$ is preserved
by $\exp$. Compute
$$\exp(\alpha_T)\exp(\alpha_{T'})\exp(\alpha_T)$$
$$=[\alpha_{(1,2)}\ldots
\alpha_{(2j-1,2j)}\alpha_{\exp(T)}].[\alpha_{\exp(T')}\alpha_{(1,2)}\ldots
\alpha_{(2i-1,2i)}].$$ $$[\alpha_{\exp(T)}\alpha_{(1,2)}\ldots \alpha_{(2j-1,2j)}]$$
$$=[\alpha_{(1,2)}\ldots\alpha_{(2j-1,2j)}].
\left(\alpha_{\exp(T)}\alpha_{\exp(T')}\alpha_{\exp(T)}\right).$$
$$[\alpha_{(2(j-i+1)-1,2(j-i+1))}\ldots\alpha_{(2j-1,2j)}][\alpha_{(1,2)}\ldots
\alpha_{(2j-1,2j)}]$$
$$=[\alpha_{(1,2)}\ldots\alpha_{(2j-1,2j)}].\alpha_{j_{\exp(T)} \exp(T')}. [\alpha_{(1,2)}\ldots
\alpha_{(2(j-i)-1,2(j-i))}]\;\;(*).$$
But $j_{\exp(T)} \exp(T')=(2(j-i+1)-1,\ldots,2j)$, so that 
$$\alpha_{j_{\exp(T)} \exp(T')}. [\alpha_{(1,2)}\ldots
\alpha_{(2(j-i)-1,2(j-i))}]=$$ $$[\alpha_{(1,2)}\ldots
\alpha_{(2(j-i)-1,2(j-i))}].\alpha_{j_{\exp(T)} \exp(T')},$$ and finally, 
$$(*)=\alpha_{(2(j-i+1)-1,2(j-i+1))}\ldots\alpha_{(2j-1,2j)}\alpha_{j_{\exp(T)}
    \exp(T')}=\exp(\alpha_{j_{\exp(T)} \exp(T')}),$$ which ends the proof of
  the first assertion of 1).\\

If now one performs, say, one simple expansion from the $i^{th}$ label,
corresponding to the expansion map $\Sigma_n\rightarrow \Sigma_{n+1}$, there
exists a lift $J_n\rightarrow J_{n+1}$: let
$\alpha=\alpha_{T_r}\ldots\alpha_{T_2}\alpha_{T_1}\in J_n$, then (supposing
$\alpha_{T_1}=\alpha_{(1,\ldots,j)}$ to simplify the notations), define first
$\exp(\alpha_{T_1})=\alpha_{(1,\ldots,j,j+1)}\alpha_{(i,i+1)}$ if $i$ belongs
to the support of $T_1$ (if not, don't modify $\alpha_{T_1}$), then define
similarly $\exp(\alpha_{T_2})$ by expanding the $\omega_{T_1}(i)^{th}$ label,
and so on. Finally you get
$\exp(\alpha):=\exp(\alpha_{T_r})\ldots\exp(\alpha_{T_1})\in J_{n+1}$, which
projects onto $\Sigma_{n+1}$ on the expansion (from the $i^{th}$ label) of the
permutation $\omega=\omega_{T_r}\ldots\omega_{T_2}\omega_{T_1}\in
\Sigma_n$. Again, it can be checked that the relations in the groups $J_n$ and
$J_{n+1}$ are preserved by this expansion map, which proves it is well-defined.\\        

2) Let $\alpha=\alpha_{T_1}\alpha_{T_2}\ldots\alpha_{T_r}\in PJ_n=Ker\, \phi$,
   $\gamma$ the combinatorial loop attached to $\alpha$, based at
   $id_n$. We claim that loop $\exp(\gamma)$ lifts to the path $(1\rightarrow
   \exp(\alpha_{T_1})\rightarrow\ldots\rightarrow
   \exp(\alpha_{T_1})\ldots\exp(\alpha_{T_r}))$, where
   $1\rightarrow\exp(\alpha_{T})$ is defined to be
$$1\rightarrow\alpha_{\exp(T)}\rightarrow\alpha_{\exp(T)}\alpha_{(1,2)}\rightarrow\ldots\rightarrow\alpha_{\exp(T)}\alpha_{(1,2)}\ldots\alpha_{(2k-1,2k)}.$$

Indeed, applying $\phi$ to this path gives precisely the loop $\exp_n(\gamma)$, as
   described in the prelimary of Theorem 4, for the embedding $\exp_n:\widetilde{{\cal
    M}}_{0,n+1}\hookrightarrow \widetilde{{\cal M}}_{0,2n+1}$. It ends at
   $\exp(\alpha_{T_1})\ldots\exp(\alpha_{T_r})=\exp(\alpha)$, so $(\exp_n)_*(\alpha)=\exp(\alpha)$.\\  

3) We use the fact that the embedding $\exp_n$ has a retraction
   $$r_n:{\widetilde{\cal M}}_{0,2n+1}(\R)\rightarrow {\widetilde{\cal
       M}}_{0,n+1}(\R),$$
which is the composition of the forgetting maps ${\widetilde{\cal
    M}}_{0,2n+1}(\R)\rightarrow {\widetilde{\cal
    M}}_{0,2n}(\R)\rightarrow\ldots\rightarrow {\widetilde{\cal
    M}}_{0,n+1}(\R)$ (cf. \cite{kn}: the maps ${\overline{\cal
    M}}_{0,n+1}(\C)\rightarrow {\overline{\cal M}}_{0,n}(\C)$ is a universal
family of $n$-pointed stable curves).\\ 
Recall also that ${\widetilde{\cal M}}_{0,n+1}(\R)$ are aspherical, so their
singular homology coincides with the homology of their fundamental groups.\\

4) The commutativity of the diagram is clear by construction. Then, injectivity of
   $PJ_n\rightarrow PJ_{2n}$ and $\Sigma_n\rightarrow \Sigma_{2n}$ implies
   injectivity of $J_n\rightarrow J_{2n}$.\\
Notice however that the retraction $(r_n)_*: PJ_{2n}\rightarrow PJ_n$ can not
   be extended to a retraction of $J_n\rightarrow J_{2n}$ (it would induce a
   retraction of $\Sigma_n\rightarrow \Sigma_{2n}$, which does not exist). $\;\;\square$

\subsection{Description of the extended groups $\tilde{\cal A}_V$ and $\tilde{\cal A}_N$}

The group $\tilde{\cal A}_V$ has a description very similar to the group $V$, by
replacing $\Sigma_n$ by the quasi-braid group $J_n$ in the definition of $V$
given at the beginning of \S 4.1. Thus, an element of $\tilde{\cal A}_V$ may be
represented by a symbol $(\alpha_1,\alpha_0, \sigma)$, where $\alpha_0$,
$\alpha_1$ are binary planar trees with $n$ leaves (for some $n$), and
$\sigma$ belongs to $J_n$. In the process of dyadic expansion for a symbol,
which may be used when composing them,
the expansion maps $\exp:J_n\rightarrow J_{n+*}$ must be used (cf. Theorem 4, {\it 1}).\\

As for the group $\tilde{\cal A}_N$, its elements are represented by symbols
$(\alpha_0,\alpha_1, q_\sigma)$ (cf. the definition of $N$) where $\sigma$
belongs to some $J_n$, and $q_\sigma$ is a collection of tree ``quasi-braid"
isomorphisms $q_i:T^{\alpha_0}_i\to T^{\alpha_1}_{\bar{\sigma}(i)}$ (where
$\bar{\sigma}\in \Sigma_n$ is the projection of $\sigma$): together with
$\sigma\in J_n$ there is a family $(\sigma_k)_{k\in \N}$ of the product
$\prod_{k\in\N} J_{2^k n}$, such that $\sigma_0=\sigma$ and $\exp(\sigma_k)$ may differ from
$\sigma_{k+1}$ by a product of some quasi-braid transpositions $\alpha_{(2i-1,2i)}$.

\subsection{Stable length, and a central extension for $N$}

Let $\alpha=\alpha_{T_1}\ldots\alpha_{T_r}$ be in the free monoid freely
generated by the generators of $J_n$. Define its length to be
$\ell_n(\alpha)=r+|T_1|+\ldots +|T_r|$, where $T_i$ is the length of the graph $G_{T_i}$.

\begin{pr}[stable length]
The length $\ell_n$ induces a well-defined group homomorphism $\ell_n:
J_n\rightarrow \Z/2\Z,$ $\ell_n(\alpha)=r+|T_1|+\ldots +|T_r|$ $mod\;
2$. Moreover, the collection $\{\ell_n,\, n\geq 1\}$ is compatible with the
direct system $\{J_n,\; \exp_n\}$, and induces a stable length $\ell_{\infty}:
J_{2^{\infty}}\rightarrow \Z/2\Z$. More generally, the length is compatible
with the dyadic expansion maps $J_{n}\rightarrow J_{n+*}$ (cf. Theorem 4, {\it
  1.}).\\
The restriction of $\ell_{\infty}$ to the infinite pure braid
group $PJ_{2^{\infty}}$ is still non-trivial. Finally, the stable length
$\ell_{\infty}$ can be extended to $\tilde{\cal A}_V$, but not to to the whole
group $\tilde{\cal A}_N$.
\end{pr}

Proof: The last two relations in the presentation of $J_n$ preserve the length
$\ell_n$. The first one ($\alpha_T ^2=1$) preserves the length $mod\; 2$
only. So $\ell_n: J_n\rightarrow \Z/2\Z$ is a well-defined group homomorphism.\\  
On the other hand, if $\alpha_T=\alpha_{(1,\ldots,k)}\in J_n$ and one performs
a simple expansion from the first leaf (to simplify the notations), then
$\exp(\alpha_T)=\alpha_{(1,\ldots,k+1)}\alpha_{(1,2)}\in J_{n+1}$, and
$\ell_n(\alpha_T)=1+k\; mod\; 2$, $\ell_{n+1}(\exp(\alpha_T))=k+1+1+2+1=\ell_n(\alpha_T)\; mod\; 2.$\\
Finally observe that the pure braid
$p=\alpha_{(1,2)}\alpha_{(2,3)}\alpha_{(1,2)}\alpha_{(1,2,3)}\in PJ_4$ %, whose
%class mod $[PJ_4,PJ_4]$ is one of the torsion-free generators of
%$H_1(PJ_4)=H_1(\widetilde{\cal M}_{0,5})$ ($\overline{\cal M}_{0,5}$ is
%the connected sum of 5 projective planes (cf. \cite{de}), hence
%$H_1(\overline{\cal M}_{0,5})=\Z^4\times \Z/2\Z$ and $H_1(\widetilde{\cal
%M}_{0,5})=\Z^4$),
 has a stable length equal to $1\;mod\; 2$.\\
Let now $g=[(\alpha_1,\alpha_0=T_{2^n},\sigma)]$ be in $\tilde{\cal A}_V$, with
$\sigma\in J_{2^n}$. Then
$\ell_{\infty}(g):=\ell_{\infty}(\bar{\sigma})$ is well-defined, since
changing of symbols for $g$ would replace $\sigma$ by some expansion of
it. On the contrary, if $g=[(\alpha_1,\alpha_0=T_{2^n},q_{\sigma})]$ is
in $\tilde{\cal A}_N$, an expansion of the symbol replaces $\sigma$ by some
expansion of it, possibly twisted by some quasi-braid transpositions
$\alpha_{(2i-1,2i)}$: but the stable length of such a quasi-braid transposition is
$1 \; mod\; 2$, and there is no possible definition for $\ell_{\infty}(g)$ (we
will confirm in the next theorem this heuristic assertion) . $\;\;\square$\\

Let now $Ker\; \ell_{\infty}$ be the kernel of the restriction of
$\ell_{\infty}$ to $PJ_{2^{\infty}}$. 

\begin{theo}[Analogue of the Euler class for $N$]
The quasi-braid extension $1\rightarrow PJ_{2^{\infty}}\rightarrow \tilde{\cal A}_N
\rightarrow N \rightarrow 1$ induces a non-trivial central extension
$$1\rightarrow \Z/2\Z\cong PJ_{2^{\infty}}/Ker\; \ell_{\infty}\longrightarrow
{\cal A}_N:=\tilde{\cal A}_N/Ker\; \ell_{\infty} \longrightarrow N\rightarrow 1,$$
which defines a non-trivial cohomology class $Eu\in H^2(N,\Z/2\Z)$.
\end{theo}

Proof: Let $g=[(\alpha_1,\alpha_0=T_{2^n},q_{\sigma})]$ be in $\tilde{\cal
  A}_N$ ($\sigma\in J_{2^n}$), and $p\in PJ_{2^{\infty}}$, represented by
  $[(\alpha_1,\alpha_1, p_1)]$, with $p_1\in PJ_{2^n}$. It follows that $g^{-1}
  p g$ is represented in $PJ_{2^n}\subset PJ_{2^{\infty}}$ by $\sigma^{-1} p_1 \sigma$, and
  $\ell_{\infty}(g^{-1} pg)=\ell_{\infty}(\sigma^{-1} p_1\sigma)=
  -\ell_{\infty}(\sigma)+\ell_{\infty}(p_1)+\ell_{\infty}(\sigma)=\ell_{\infty}(p)$,
  using Proposition 5.
  This proves in particular that $[\tilde{\cal A}_N, PJ_{2^{\infty}}]\subset Ker\;
  \ell_{\infty}$: so, $Ker\; \ell_{\infty}$ is normal in $\tilde{\cal A}_N$, and the extension is central.\\
Suppose the extension is trivial: then the embedding $i:\Z/2\Z\rightarrow
  \tilde{\cal A}_N$ would admit a retraction $r$. But then, the composition
  $PJ_{2^{\infty}}\rightarrow \tilde{\cal A}_N \rightarrow \tilde{\cal A}_N/Ker\;
  \ell_{\infty}\stackrel{r}{\rightarrow} PJ_{2^{\infty}}/Ker\;
  \ell_{\infty}\cong\Z/2\Z$ would be $\ell_{\infty}$, proving that the
  composition of the last two morphisms is a prolongation to $\tilde{\cal A}_N$
  of the stable length morphism, which is impossible, admitting the assertion of the
  preceding proposition.\\
However, we now give another proof of the non-triviality independent of the
  assertion on the stable length: the method consists in writing the generator
  of the kernel $\Z/2\Z\cong PJ_{2^{\infty}}/Ker\; \ell_{\infty}$ as a product
  of commutators in $\tilde{\cal A}_N/Ker\; \ell_{\infty}$, i. e. finding a
  pure quasi-braid with length $1\;mod\; 2$ as a product of commutators in
  $\tilde{\cal A}_N$ (which definitely proves that $\ell_{\infty}$ can not be
  extended to $\tilde{\cal A}_N$).\\
Let $\tau\in V\subset N$ be the transposition defined by the symbol  \makebox{\begin{picture}(0,0)%
\includegraphics{tau.pstex}%
\end{picture}%
\setlength{\unitlength}{1906sp}%
\begingroup\makeatletter\ifx\SetFigFont\undefined%
\gdef\SetFigFont#1#2#3#4#5{%
  \reset@font\fontsize{#1}{#2pt}%
  \fontfamily{#3}\fontseries{#4}\fontshape{#5}%
  \selectfont}%
\fi\endgroup%
\begin{picture}(1497,1374)(2566,-1423)
\put(2566,-1411){\makebox(0,0)[lb]{\smash{\SetFigFont{6}{7.2}{\rmdefault}{\mddefault}{\updefault}$a$}}}
\put(3646,-1411){\makebox(0,0)[lb]{\smash{\SetFigFont{6}{7.2}{\rmdefault}{\mddefault}{\updefault}$b$}}}
\end{picture}}, the leaves $a$ and
  $b$ being permuted.\\
Let $\alpha\in Aut(T_2)\subset N$ be defined by the symbol \makebox{\begin{picture}(0,0)%
\includegraphics{alpha.pstex}%
\end{picture}%
\setlength{\unitlength}{1906sp}%
\begingroup\makeatletter\ifx\SetFigFont\undefined%
\gdef\SetFigFont#1#2#3#4#5{%
  \reset@font\fontsize{#1}{#2pt}%
  \fontfamily{#3}\fontseries{#4}\fontshape{#5}%
  \selectfont}%
\fi\endgroup%
\begin{picture}(2352,2487)(1711,-2536)
\put(1711,-2536){\makebox(0,0)[lb]{\smash{\SetFigFont{6}{7.2}{\rmdefault}{\mddefault}{\updefault}...}}}
\put(2431,-1411){\makebox(0,0)[lb]{\smash{\SetFigFont{6}{7.2}{\rmdefault}{\mddefault}{\updefault}$a$}}}
\put(3646,-1411){\makebox(0,0)[lb]{\smash{\SetFigFont{6}{7.2}{\rmdefault}{\mddefault}{\updefault}$b$}}}
\put(3196,-1861){\makebox(0,0)[lb]{\smash{\SetFigFont{6}{7.2}{\rmdefault}{\mddefault}{\updefault}2}}}
\put(1981,-1861){\makebox(0,0)[lb]{\smash{\SetFigFont{6}{7.2}{\rmdefault}{\mddefault}{\updefault}  1}}}
\end{picture}}  (the permutations of the leaves indicated by the arrows must
be read {\it from bottom to top}). Set $\gamma:=\tau\alpha$. Then $\gamma=$ \makebox{\begin{picture}(0,0)%
\includegraphics{gamma.pstex}%
\end{picture}%
\setlength{\unitlength}{1906sp}%
\begingroup\makeatletter\ifx\SetFigFont\undefined%
\gdef\SetFigFont#1#2#3#4#5{%
  \reset@font\fontsize{#1}{#2pt}%
  \fontfamily{#3}\fontseries{#4}\fontshape{#5}%
  \selectfont}%
\fi\endgroup%
\begin{picture}(2352,2487)(1711,-2536)
\put(1711,-2536){\makebox(0,0)[lb]{\smash{\SetFigFont{6}{7.2}{\rmdefault}{\mddefault}{\updefault}...}}}
\put(2431,-1411){\makebox(0,0)[lb]{\smash{\SetFigFont{6}{7.2}{\rmdefault}{\mddefault}{\updefault}$a$}}}
\put(3646,-1411){\makebox(0,0)[lb]{\smash{\SetFigFont{6}{7.2}{\rmdefault}{\mddefault}{\updefault}$b$}}}
\put(3196,-1861){\makebox(0,0)[lb]{\smash{\SetFigFont{6}{7.2}{\rmdefault}{\mddefault}{\updefault}2}}}
\put(1981,-1861){\makebox(0,0)[lb]{\smash{\SetFigFont{6}{7.2}{\rmdefault}{\mddefault}{\updefault}  1}}}
\end{picture}}, and it appears that $\gamma$ and $\alpha$ are conjugated by the ``translation"
  $\delta=$ \makebox{\begin{picture}(0,0)%
\includegraphics{delta.pstex}%
\end{picture}%
\setlength{\unitlength}{1906sp}%
\begingroup\makeatletter\ifx\SetFigFont\undefined%
\gdef\SetFigFont#1#2#3#4#5{%
  \reset@font\fontsize{#1}{#2pt}%
  \fontfamily{#3}\fontseries{#4}\fontshape{#5}%
  \selectfont}%
\fi\endgroup%
\begin{picture}(4842,1591)(975,-1644)
\put(3601,-781){\makebox(0,0)[lb]{\smash{\SetFigFont{6}{7.2}{\rmdefault}{\mddefault}{\updefault},}}}
\end{picture}
}. Precisely, we have $\gamma=\delta\alpha\delta^{-1}$, or equivalently,
  $\tau=[\delta,\alpha]$.\\
We now lift $\tau$, $\delta$ and $\alpha$ in $\tilde{\cal A}_N$ in an obvious
  way: $\tau$ is lifted in $\tilde{\tau}$ (same symbol coupled with
  $\alpha_{(12)}\in J_{2^{\infty}}$), $\delta$ in $\tilde{\delta}$ (same symbol coupled with
  $1\in J_{2^{\infty}}$), and $\alpha$ lifted in $\tilde{\alpha}$ (same symbol
  coupled with the sequence $\alpha_0=1$, $\alpha_1=1$,
  $\alpha_2=\alpha_{(12)}$,  $\alpha_{k+1}=\exp(\alpha_k)\alpha_{(12)}$, $k\in \N$): \makebox{\begin{picture}(0,0)%
\includegraphics{alpha2.pstex}%
\end{picture}%
\setlength{\unitlength}{1906sp}%
\begingroup\makeatletter\ifx\SetFigFont\undefined%
\gdef\SetFigFont#1#2#3#4#5{%
  \reset@font\fontsize{#1}{#2pt}%
  \fontfamily{#3}\fontseries{#4}\fontshape{#5}%
  \selectfont}%
\fi\endgroup%
\begin{picture}(2802,2847)(1339,-2896)
\put(2431,-1411){\makebox(0,0)[lb]{\smash{\SetFigFont{6}{7.2}{\rmdefault}{\mddefault}{\updefault}$a$}}}
\put(3646,-1411){\makebox(0,0)[lb]{\smash{\SetFigFont{6}{7.2}{\rmdefault}{\mddefault}{\updefault}$b$}}}
\put(3196,-1861){\makebox(0,0)[lb]{\smash{\SetFigFont{6}{7.2}{\rmdefault}{\mddefault}{\updefault}2}}}
\put(1981,-1861){\makebox(0,0)[lb]{\smash{\SetFigFont{6}{7.2}{\rmdefault}{\mddefault}{\updefault}  1}}}
\put(4141,-961){\makebox(0,0)[lb]{\smash{\SetFigFont{6}{7.2}{\rmdefault}{\mddefault}{\updefault}$\alpha_0=1\in J_2$}}}
\put(4141,-1411){\makebox(0,0)[lb]{\smash{\SetFigFont{6}{7.2}{\rmdefault}{\mddefault}{\updefault}$\alpha_1=1\in J_4$}}}
\put(4141,-1861){\makebox(0,0)[lb]{\smash{\SetFigFont{6}{7.2}{\rmdefault}{\mddefault}{\updefault}$\alpha_2=\alpha_{(12)}\in J_8$}}}
\put(4141,-2266){\makebox(0,0)[lb]{\smash{\SetFigFont{6}{7.2}{\rmdefault}{\mddefault}{\updefault}$\alpha_3=\exp(\alpha_{(12)})\alpha_{(12)}\in J_{16}$}}}
\put(1351,-2896){\makebox(0,0)[lb]{\smash{\SetFigFont{6}{7.2}{\rmdefault}{\mddefault}{\updefault}...}}}
\put(4141,-2761){\makebox(0,0)[lb]{\smash{\SetFigFont{6}{7.2}{\rmdefault}{\mddefault}{\updefault}$\alpha_{k+1}=\exp(\alpha_k)\alpha_{(12)}\in J_{2^{k+2}}$}}}
\end{picture}}\\
 Clearly, the same relation as in $N$ holds:
 $\tilde{\tau}=[\tilde{\delta},\tilde{\alpha}]\in \tilde{\cal A}_N$.\\
On the other hand, $\tau=$ \makebox{\begin{picture}(0,0)%
\includegraphics{tau2.pstex}%
\end{picture}%
\setlength{\unitlength}{1906sp}%
\begingroup\makeatletter\ifx\SetFigFont\undefined%
\gdef\SetFigFont#1#2#3#4#5{%
  \reset@font\fontsize{#1}{#2pt}%
  \fontfamily{#3}\fontseries{#4}\fontshape{#5}%
  \selectfont}%
\fi\endgroup%
\begin{picture}(1775,2318)(2288,-2367)
\put(2431,-1411){\makebox(0,0)[lb]{\smash{\SetFigFont{6}{7.2}{\rmdefault}{\mddefault}{\updefault}$a$}}}
\put(3646,-1411){\makebox(0,0)[lb]{\smash{\SetFigFont{6}{7.2}{\rmdefault}{\mddefault}{\updefault}$b$}}}
\put(2341,-2041){\makebox(0,0)[lb]{\smash{\SetFigFont{6}{7.2}{\rmdefault}{\mddefault}{\updefault}1}}}
\put(2971,-2041){\makebox(0,0)[lb]{\smash{\SetFigFont{6}{7.2}{\rmdefault}{\mddefault}{\updefault}2}}}
\put(3331,-2041){\makebox(0,0)[lb]{\smash{\SetFigFont{6}{7.2}{\rmdefault}{\mddefault}{\updefault}3}}}
\put(3871,-2041){\makebox(0,0)[lb]{\smash{\SetFigFont{6}{7.2}{\rmdefault}{\mddefault}{\updefault}4}}}
\end{picture}}  may also be written as the product
$\tau=\tau_1\tau_2$, where $\tau_1$ exchanges the leaves 1 and 3 (keeping 2
and 4 fixed) and $\tau_2$ exchanges the leaves 2 and 4 (keeping 1
and 3 fixed). We note abusively $\tau_1=(13)$, $\tau_2=(24)$. Introducing
$\sigma\in V$
defined by $\sigma=(12)(34)$, we have $\tau_2=\sigma\tau_1\sigma$, and
$\tau=[\tau_1,\sigma]$.\\
We then lift $\tau_1$ and $\sigma$ in $\tilde{\cal A}_N$ by
$\tilde{\tau_1}=\alpha_{(123)}$, $\tilde{\sigma}= \alpha_{(12)}\alpha_{(34)}$.
Now $[\tilde{\tau_1},\tilde{\sigma}]$ differs from
$\tilde{\tau}=[\tilde{\delta},\tilde{\alpha}]$ by a pure quasi-braid
$$p=[\tilde{\tau_1},\tilde{\sigma}][\delta,\alpha]=[\alpha_{(123)},\alpha_{(12)}\alpha_{(34)}]\exp(\alpha_{(12)})$$
$$=[\alpha_{(123)},\alpha_{(12)}\alpha_{(34)}]\alpha_{(12)}\alpha_{(34)}\alpha_{(1234)}=\alpha_{(123)}\alpha_{(12)}\alpha_{(34)}\alpha_{(123)}\alpha_{(1234)}.$$
The miracle is that $\ell_{\infty}(p)=1\;mod\;2$ as desired. $\;\;\square$\\

\begin{co} The 2-cycle $\omega$ defined by the relation
  $[\tau_1,\sigma][\alpha,\delta]=1\in N$ is non-trivial and verifies
  $(Eu,[\omega])=1$, where $Eu\in H^2(N,\F_2)$ is the cohomology class of the
extension of $N$.
\end{co}
This is an immediate application of the following
\begin{lem}
Let $G$ be a perfect group, $A\rightarrow\hat{G}\rightarrow G$ a central
extension of $G$ with kernel an abelian group $A$,
$c\in H^2(G,A)$ the associated cohomology class. If $\omega$ is a 2-cycle of
$G$ associated to a relation $1=\prod_i [g_i,h_i]$ in $G$, then $(c,[\omega])=a\in A$,
where $a$ is computed as $a=\prod_i [\hat{g_i},\hat{h}_i]$, for any choices of
lifts $\hat{g_i}$, $\hat{h}_i$ of $g_i$, $h_i$.
\end{lem}

The proof is easy by describing the extension $A\rightarrow\hat{G}\rightarrow
G$ as a push-out of the universal central extension $H_2(G)\rightarrow
G^{univ}\rightarrow G$ as defined by Hopf Theorem. Note that the map $H_2(G)\rightarrow A$ defining the
push-out corresponds to the cohomology class $c$ under the isomorphism
$H^2(G,A)\cong Hom(H_2(G),A)$.\\

So, in our case, $(Eu,[\omega])=\ell_{\infty}(p)=1\;mod\; 2$. $\;\;\square$\\

\subsection{Euler extension for the model ${\cal N}$ of Neretin group}  

 We may now prefer considering the analogous construction for Neretin group ${\cal
  N}$ related to the unrooted regular dyadic tree ${\cal T}_2$, since it
  naturally contains the tree automorphism group $Aut({\cal T}_2)$ and the
  group $PSL(2,\Q_2)$. 

\begin{lem}
\begin{enumerate}
\item The subgroup $Ker\, \ell_n$ of $PJ_n$ is normal in $Q_{n+1}$, and the central
extension $0\rightarrow \Z/2\Z=PJ_n/Ker\, \ell_n\longrightarrow Q_{n+1}/Ker\,
\ell_n\rightarrow \Z/2\Z\rightarrow 0$ is trivial. Consequently, the length
$\ell_n$ may be extended to a length $\bar{\ell}_{n+1}$ on $Q_{n+1}$.
\item Let $\overline{\exp}_n:Q_n\rightarrow Q_{2n}$ be the morphism induced by the
embedding $\overline{\cal M}_{0,n}(\R)\\\hookrightarrow \overline{\cal
  M}_{0,2n}(\R)$. Then $\overline{\exp}_n$ maps $\hat{a}\in Q_{n}$ to a
pure quasi-braid in $PJ_{2n-1}$. Consequently, the group $Q_{3.2^{\infty}}$ is isomorphic to
$PJ_{2^{\infty}}$, and there is a well-defined stable length $\bar{\ell}_{\infty}$ on $Q_{3.2^{\infty}}$.
\end{enumerate}
\end{lem}

Proof: 1. As we have observed, there is a pure quasi-braid
$p=\hat{\alpha}j_S(\hat{\alpha})$ such that $\hat{a}^{-1}=p^{-1}\hat{a}$. For
each $q\in PJ_n$, $\hat{a}^{-1}q\hat{a}$ is the pure quasi-braid equal to
$\hat{\alpha} j_S( p^{-1}) j_S( q) j_S(\hat{\alpha})$. Since $\ell_n(p)=0$, we get
$\ell_n(\hat{a}^{-1}q\hat{a})=\ell_n(q)$, and the central extension is
well-defined. It is trivial, since $\hat{a}^2=p$ belongs to $Ker\,
\ell_n$. Consequently, we may set
$\bar{\ell}_{n+1}(\hat{a})=0$, this extends $\ell_n$, as a
morphism, from
$PJ_n$ to $Q_{n+1}$.\\

2. Let us choose
 $\hat{\alpha}=\alpha_{(1,2)}.(\alpha_{(2,3)}\alpha_{(1,2)}).(\alpha_{(3,4)}\alpha_{(2,3)}\alpha_{(1,2)}).\;\;
 \ldots\;\; 
 .(\alpha_{(n-1,n)}\\\alpha_{(n-2,n-1)}\ldots \alpha_{(1,2)})\in J_n$, to define
 $\hat{a}\in Q_{n+1}$. If we proceed to a simple expansion from the last
 puncture, we get an embedding $\overline{\cal M}_{0,n+1}(\R)\hookrightarrow
 \overline{\cal  M}_{0,n+2}(\R)$ and an induced morphism $Q_{n+1}\rightarrow
 Q_{n+2}$. It is easy to see that this morphism maps $\hat{a}$ on the pure
 quasi-braid $\hat{\alpha}\alpha_{(1,\ldots,n)}\in PJ_{n+1}$. This of course
 generalizes to an arbitrary number of expansions. Consequently, the stable
 length $\bar{\ell}_{\infty}$ is well-defined on $Q_{3.2^{\infty}}$, since
 each element of $Q_{3.2^{\infty}}$ may be represented by a pure quasi-braid
 (we just need the definition of the
 $\ell_n$'s on the pure-braid groups, not the $\bar{\ell}_n$'s).
 $\;\;\square$

\begin{de}[Euler class of ${\cal N}$] 

The extension $Q_{3.2^{\infty}}\rightarrow \bar{\cal A}_{\cal
  N}\rightarrow {\cal N}$ induces the central extension $$1\rightarrow
  \Z/2\Z\cong Q_{3.2^{\infty}}/Ker\, \bar{\ell}_{\infty}\longrightarrow
  \bar{\cal A}_{\cal N} / Ker\; \bar{\ell}_{\infty}\longrightarrow {\cal N} \rightarrow 1$$
which we call the Euler class of ${\cal N}$.
\end{de}
   
\subsection{Euler cocycle}

Let ${\cal R}$ be the ring of $\Z/2\Z$-valued sequences, divided by the ideal
of almost zero sequences: ${\cal R}= (\Z/2\Z)^{\N}/(\Z/2\Z)^{(\N)}$. Denote by
$1_{\cal R}$ its unit.\\
For each ${\tilde f}$ in $\tilde{\cal A}_N$ defined by a symbol of the form
$(\alpha_1,\alpha_0=T_{2^n}, q_{\sigma})$ (cf. \S 3.2), there is a family
$(\sigma_k)_{k\geq n}$, $\sigma_k \in J_{2^k}$, with $\sigma_n=\sigma$ and
$\sigma_{k+1}$ differing from $\sigma_k$ by a product of quasi-braid
transpositions. So there is a well-defined function
$$\tilde{\ell}:\tilde{\cal A}_N\rightarrow {\cal R},\;\;\; \tilde{f}\mapsto \tilde{\ell}(\tilde{f}),$$
where $\tilde{\ell}(\tilde{f})$ is represented by the sequence
$\tilde{\ell}(\tilde{f})_k=\ell_{\infty}(\sigma_k)\in \Z/2\Z$ for $k\geq n$
and, say, $\tilde{\ell}(\tilde{f})_k=0$ for $k=0,\ldots, n-1$.\\

Denote by $j:\Z/2\Z=PJ_{2^{\infty}}/Ker\, \ell_{\infty}\hookrightarrow {\cal
  R}$ the natural embedding.

\begin{lem}
For all ${\tilde f}$ in $\tilde{\cal A}_N$ and $p\in PJ_{2^{\infty}}$,
$\tilde{\ell}(p \tilde{f})=\tilde{\ell}(p)+
\tilde{\ell}(\tilde{f})=\tilde{\ell}(\tilde{f}p)$, and $\tilde{\ell}$
induces a function $\ell: {\cal A}_N=\tilde{\cal A}_N/Ker\, \ell_{\infty}\rightarrow {\cal R}$ such that for all $f\in {\cal A}_N$,  
$$\ell(T.f)=1_{\cal R} + \ell(f),$$   
where $T=j(1)$.
\end{lem}

Proof: Choose $n$ such that $p$ may be represented in $PJ_{2^n}$, and ${\tilde
  f}$ represented by a symbol $(T_1,T_{2^n},q_{\sigma})$. Write
  $p=[(T_{2^n},T_{2^n},\alpha)]=[(T_1,T_1,\beta)]$, where $\alpha$ and $\beta$
  possess a common expansion in some $PJ_{2^m}$, $m\geq n$. Then ${\tilde
  f}p=[(T_1,T_{2^n},q_{\tau})]$ with $\tau=\tau_n=\sigma\alpha$,
  $\tau_k=\sigma_k \exp_k(\alpha)$ for $k\geq n$, $p{\tilde f}=[(T_1,T_{2^n},q_{\upsilon})]$
  with $\upsilon=\upsilon_n=\beta\sigma$,
  $\upsilon_k=\exp_k(\beta)\sigma_k$. Since
  $\ell_{\infty}(\exp_k(\alpha))=\ell_{\infty}(\exp_k(\beta))=\ell_{\infty}(p)$
  $\forall k\geq n$,
  the proof is done. Note that
  $\tilde{\ell}(p)=j(\ell_{\infty}(p))$. $\;\;\square$\\

Since $N$ is perfect, the injection $j:\Z/2\Z\hookrightarrow {\cal   R}$
induces an injective morphism $j_*: H^2(N,\Z/2\Z) \hookrightarrow H^2(N,{\cal
  R})$.

\begin{theo}
The image by $j_*$ of the Euler class $Eu\in H^2(N,\Z/2\Z)$ is the cohomology
class  of the well-defined cocycle $c:N\times N\rightarrow {\cal R}$ defined
by
$$c(f,g)=\ell(\bar{f}\bar{g})-\ell(\bar{f})-\ell(\bar{g})$$
where $\bar{f}$, $\bar{g}$ are any lifts to ${\cal A}_N$ of $f$ and $g$
respectively.
\end{theo}

Proof: First the fact that the cocycle $c$ is well-defined follows from the
equivariant relation of Lemma 3.\\
Let $\omega$ be a 2-cycle of $N$. It is associated to a relation
$\prod_{i=1}^p  [f_i,g_i]=1\in N$, and may be written
$$\omega=\sum_{i=1}^p (f_i,g_i)-(g_i,f_i)-(g_i f_i,(g_i f_i)^{-1})+(f_i
g_i,(g_i f_i)^{-1})$$
$$+\sum_{i=1}^{p-1} ([f_1,g_1]\ldots [f_i,g_i],
[f_{i+1},g_{i+1}]),$$
and it follows that $([c],[\omega])=\ell(\prod_{i=1}^p 
[\bar{f}_i,\bar{g}_i])$. Now $\prod_{i=1}^p [\bar{f}_i,\bar{g}_i]=\alpha\, mod\,
Ker\, \ell_{\infty}$, for some $\alpha\in PJ_{2^{\infty}}$, and
$\ell(\prod_{i=1}^p
[\bar{f}_i,\bar{g}_i])=\tilde{\ell}(\alpha)=j(\ell_{\infty}(\alpha))$. But as
mentioned in \S 3.3, Lemma 1, $\ell_{\infty}(\alpha)=(Eu,[\omega])$, so that
$([c],[\omega])=j((Eu,[\omega]))=(j_*Eu,[\omega])$, which proves $[c]=j_*Eu$,
since $H^2(N,{\cal R})=Hom(H_2(N), {\cal R})$. $\;\;\square$

\subsection{The analogy with the Euler class of homeomorphism groups of the
  circle}

1. First (but a little naive) evidence to think of this new class on ${\cal N}$ as the analogue of the Euler
class of Thompson group ${\cal T}$ (cf. \cite{gh-se}):\\
The latter is obtained by lifting to
$\R$ the action of ${\cal T}$ on the circle. It is simply the restriction to
${\cal T}$ of the Euler class of the group $H\widetilde{omeo}^+(S^1)$ of
orientation-preserving homeomorphisms of the circle, namely the class of the extension
$0\rightarrow \Z \rightarrow H\widetilde{omeo}^+(S^1)\rightarrow Homeo^+(S^1)
\rightarrow 1,$ where $H\widetilde{omeo}^+(S^1)$, the universal cover of the
group $Homeo^+(S^1)$, is the group of $\Z$-equivariant homeomorphisms of $\R$
(viewing $S^1$ as $\R/\Z$). At first
sight, one should be temptated to accept the boundary $\partial {\cal T}_2$ of the dyadic infinite
tree, on which ${\cal N}$ continuously acts, as the $p$-adic analogue of the circle
(since $\partial {\cal T}_2\cong \Q_2 P^1$ and $S^1\cong \R P^1$). But $\Q_2 P^1$
being totally disconnected, this kills any
hope to do some topology. At this point, the modular tower reveals to us as
the appropriate space related to ${\cal N}$, and the non-triviality of its homotopy type
generates a non-trivial cohomology class for ${\cal N}$ -- just as the homotopy type of the
circle generates Thompson group ${\cal T}$'s Euler class.\\
However it could be objected (and it was objected by V. Sergiescu to the
author) that the discrete Godbillon-Vey type class $gv$ of
Thompson group ${\cal T}$ (cf. \cite{gh-se}) may also be derived from a topological
extension of ${\cal T}$, namely the Greenberg-Sergiescu braid
extension $1\rightarrow B_{\infty}\longrightarrow A\longrightarrow{\cal
  T}\rightarrow 1$
(cf. \cite{gr-se}),
by abelianization of the kernel. So we need a more convincing argument.\\

2. The more serious evidence relies on the relation between the Euler class
   of $Homeo^+(S^1)$ and {\it bounded cohomology}.\\
Indeed, recall from \cite{ba-gh} that there is a cocycle $ord$ (``cyclic order") on $Homeo^+(S^1)$ whose class is twice the class of $0\rightarrow \Z \rightarrow H\widetilde{omeo}^+(S^1)\rightarrow Homeo^+(S^1)
\rightarrow 1,$ such that
   $ord(f,g)=\Phi_{ord}(\tilde{f}\circ\tilde{g})-\Phi_{ord}(\tilde{f})-\Phi_{ord}(\tilde{g})$,
   for any lifts $\tilde{f}$, $\tilde{g}$ in $H\widetilde{omeo}^+(S^1)$ of $f$
   and $g$ respectively, and $\Phi_{ord}(\tilde{f})=2E(\tilde{f}(0))\in \Z $,
   where $E$ is a certain modified integer part function. And the embedding of
   coefficients $\Z\rightarrow\R$ maps the integral Euler class to the class
   of the real cocycle
   $eu(f,g)=\tau(\tilde{f}\circ\tilde{g})-\tau(\tilde{f})-\tau(\tilde{g})$,
   where $\tau(\tilde{f})=\lim_{n\to +\infty} \frac{\tilde{f}^n (0)}{n}$
   is the {\it translation number of Poincar\'e}. The cocycles $ord$ and $eu$
   are induced by the boundary of an unbounded function on
   $H\widetilde{omeo}^+(S^1)$, and the class of $eu$ stands in the bounded
   cohomology group $H^2_b (Homeo^+(S^1),\R)$. We think the analogy with
   our class $Eu$ is very suggestive, replacing $\Z$, $\R$ and $\tau$ (or $\Phi_{ord}$)
   by $\Z/2\Z$, ${\cal R}$ and $\ell$ respectively.\\

%Proof: If it were trivial, $\bar{\ell}_{\infty}$ could be extendable to a
%group morphism on $\overline{Aut}({\cal T}_2)$. We now show it is impo%ssible.
%Consider the element $\alpha$ in $\overline{Aut}(T_2)\subset\overline{Aut}({\cal T}_2)$
%($T_2$ being the rooted planar infinite tree, union of the finite rooted trees
%$T_{2^n}$, embedded in the unrooted tree ${\cal T}_2$) defined by $\alpha=(\alpha_k)_k\in \prod_{k=0}^{+\infty} J_{2^k}$,
%$\alpha_{k+1}=\exp(\alpha_k)\alpha_{(1,2)}$ if $k\geq 0$, $\alpha_0=1$. Note
%by the way that the extension over ${Aut}(T_2)$ is trivial, so that $\overline{Aut}(T_2)$ is simply isomorphic to ${Aut}(T_2)\times PJ_{2^{\infty}}$.\\
%Let
%$\beta=(\beta_k)_k\in J_{2^{\infty}}$ be defined by $\beta_0=\alpha_{(1,2)}$
%and $\beta_{k+1}=\exp(\beta_k)$. Then $\gamma:=\beta\alpha$ is such that
%$\gamma_k=\alpha_{k+1}$, $k\geq 0$. It follows that $\gamma$ and $\alpha$ are
%conjugated in $\overline{Aut}({\cal T}_2)$, so
%$\bar{\ell}_{\infty}(\alpha)=\bar{\ell}_{\infty}(\gamma)$. But
%$\bar{\ell}_{\infty}(\gamma)=\bar{\ell}_{\infty}(\alpha)+\bar{\ell}_{\%infty}(\beta)$,
%with $\bar{\ell}_{\infty}(\beta)=1$, hence the contradiction. PREUVE FAUSSE!
%EN EFFET, JE SOUS-ENTENDS QUE LE PROLONGEMENT, S'IL EXISTE, EST COMPATIBLE
%AVEC LE PROLONGEMENT EVIDENT A $J_{2^{\infty}}$ $\;\;\square$\\ 

3. In \cite{kap2} we have introduced a so-called analogue of the Virasoro
extension of $Diff^+(S^1)$, the orientation-preserving diffeomorphism group of
the circle, for the discrete group ${\cal N}$. It is defined as follows:
Let $PAut({\cal T}_2)$ be the group of bijections on the set of vertices
${\cal T}_2^0$ of
the regular dyadic tree ${\cal T}_2$, which induce a simplicial action outside some finite subtree of
${\cal T}_2$. Each bijection in $PAut({\cal T}_2)$ induces an element of the
group ${\cal N}$ by
forgetting the action on finite subtrees, and there is a short exact sequence
$$1\rightarrow \Sigma_{\infty}\longrightarrow PAut({\cal T}_2)\longrightarrow
{\cal N}\rightarrow 1,$$ where $\Sigma_{\infty}$ is the group of finitely supported
permutations on the set ${\cal T}_2^0$. Dividing by the alternating group ${\mathfrak
  A}_{\infty}$ now provides the central extension
$$1\rightarrow \Z/2\Z\cong \Sigma_{\infty}/{\mathfrak
  A}_{\infty}\longrightarrow PAut({\cal T}_2)/{\mathfrak A}_{\infty}\longrightarrow
{\cal N}\rightarrow 1,$$ 
which is non-trivial (cf. \cite{kap2}). Denote by $Gv$ the cohomology class of
this extension. It is certainly different from the Euler class $Eu$ we have
just defined (a rigorous proof would require the construction of a
non-trivial 2-cycle
which gives different values once evaluated on $Eu$ and $Gv$, but we are not
able to find it). Though the analogy with
   Thompson group ${\cal T}$ is striking -- indeed, $H^2({\cal T},\Z)=\Z
   gv\oplus \Z eu$, cf. \cite{gh-se}, there is no natural relation between the
   cohomology classes on Neretin group and Thompson group, since the embedding
   ${\cal T}\hookrightarrow {\cal N}$ factors through Thompson group ${\cal
   V}$ related to the Cantor set, which has no cohomology in degree 2.\\    

\section{Further results}
 In this section we give further results on the groups concerned with this
 article, with sketch of proofs only, to avoid to lengthen the material too
 much.\\
  
\noindent{\bf 4.1 Vanishing of the Euler class on the subgroup
 $PGL(2,\Q_2)$}. We can not avoid the natural question of understanding the Euler class on the
subgroup $PGL(2,\Q_2)$ of the automorphism group of the tree. Though the
non-central extension of $PGL(2,\Q_2)$ by the pure quasi-braid group does not
seem to be trivial, we prove that the induced central extension is, without
however exhibiting an ``obvious" section or retraction. In particular, our
Euler class is not related with the Euler cocycle of Barge constructed on  $PSL(2,k)$, for every field $k$, with
values in the Witt group $W(k)$ (cf. \cite{ba}). This negative result means that the nature of
our class is not arithmetic.

\begin{theo}
The Euler class restricted to the subgroup $PGL(2,\Q_2)$ vanishes.
\end{theo}

Proof (sketch): The proof first exploits the knowledge of
$H_2(GL(2,\Q_2),\Z)=K_2(\Q_2)\oplus H_2(\Q_2^*,\Z)$, where $K_2(\Q_2)$ is Milnor's $K_2$, and the fact that the
pull-back of $Eu$ on $GL(2,\Q_2)$ gives zero once evaluated on
$H_2(GL(2,\Q_2),\Z)$; second, it relies on the
evaluation of the restriction of the class $Eu$ on $PGL(2,\Q_2)$ on a 2-cycle
associated to a certain way of writing $-I_2\in GL(2,\Q_2)$ as a product of
two commutators. The result of the evaluation, which requires a very heavy
computation, is zero again; finally we check that $Eu$ does not come from a
class of $Ext(H_1(PGL(2,\Q_2)),\Z/2\Z)$ either. $\;\;\square$\\

\noindent{\bf 4.2. Topological monoids related to the short exact sequence $1\rightarrow PJ_{2^{\infty}}\longrightarrow J_{2^{\infty}} \longrightarrow
\Sigma_{2^{\infty}}\rightarrow 1$}. The aim of this section is essentially to
mention that the groups $PJ_n$ and $J_n$ provide after stabilisation new
examples of groups, namely $PJ_{2^{\infty}}$ and $J_{2^{\infty}}$, which are
    homologically equivalent to loop spaces. Recall from \cite{mc-se} that the
    group $\Sigma_{\infty}$ of finitely supported permutation on a countable set is homologically equivalent to one connected component $Q$ of
 the infinite loop space $\Omega^{\infty}S^{\infty}$
 (Barratt-Priddy-Quillen theorem), whereas $\Sigma_{2^{\infty}}$ is
 homologically equivalent to the localisation $Q[2^{-1}]$, because the expansion map
 $B\Sigma_{2^k}\rightarrow B\Sigma_{2^{k+1}}$ (at the classifying space level)
 corresponds to multiplying by $2$ in the $H$-space structure of $Q$.   
   We first recall the construction of a telescope related to the group
   $\Sigma_{2^{\infty}}$, inspired from \cite{mc-se}:\\
Consider the group morphisms
$\Sigma_{2^k}\times\Sigma_{2^l}\rightarrow \Sigma_{2^{k+l}}$,
$(\sigma,\tau)\mapsto \sigma(\tau,\ldots,\tau)$, where
$\sigma(\tau,\ldots,\tau)$ denotes the wreath-product of $\sigma$ with $2^k$
copies of $\tau$. It may be seen as a restriction of the map (not a
morphism) $\Sigma_r\times\Sigma_{n_1}\times\ldots\times\Sigma_{n_r}\rightarrow
\Sigma_{n_1+\ldots+n_r}$, $(\sigma,\tau_1,\ldots,\tau_r)\mapsto
\sigma(\tau_1,\ldots,\tau_r)$ (with $r=2^k$, $n_1=\ldots=n_r=2^l$).\\
They provide the disjoint union $M_{\Sigma}:=\coprod_{k\in \N} B\Sigma_{2^k}$ with a
structure of topological monoid. Here, $B$ denotes the classifying space
functor. The unit of the monoid $\pi_0(M_{\Sigma})\cong \N$ of connected
components of $M_{\Sigma}$ comes from the identity $id_2$ in $\Sigma_2$. Now
form the telescope
$$M_{\Sigma}\stackrel{.id_2}{\rightarrow}M_{\Sigma}\stackrel{.id_2}{\rightarrow}\ldots,$$
where $.id_2$ is right multiplication induced by $id_2$. The crucial point is that the
right multiplication by $id_2$ is equivalent to the dyadic expansion
process. So, the inductive limit of the telescope, $(M_{\Sigma})_{\infty}$, is
isomorphic to $\Z\times B\Sigma_{2^{\infty}}$.\\

 We now lift the preceding construction to the
quasi-braid groups.
\begin{pr}
The group morphisms $\Sigma_{2^k}\times\Sigma_{2^l}\rightarrow
\Sigma_{2^{k+l}}$ can be lifted to the quasi-braid groups as morphisms
$$J_{2^k}\times J_{2^l}\rightarrow
J_{2^{k+l}},$$ and there is a corresponding topological monoid
$M_{J}:=\coprod_{k\in \N} BJ_{2^k}$, as well as a telescope $(M_{J})_{\infty}\cong \Z\times
BJ_{2^{\infty}}$. 
\end{pr}

Proof: We define the morphism $J_{2^k}\times J_{2^l}\rightarrow J_{2^{k+l}},$
by 2 commuting morphisms $J_{2^k}\rightarrow J_{2^{k+l}}$ and
$J_{2^l}\rightarrow J_{2^{k+l}}$. The first one is the dyadic expansion
morphism of Theorem 3
(iterated $l$ times), $\exp^{l}$; as for the second one, let $\alpha_T$ be a generator
of $J_{2^l}$, where $T$ is thought of as a $2^l$-labeled tree with a unique
internal edge. Let $T_{2^k}$ be the star with $2^k$ leaves. Graft the tree $T$
to the $i^{th}$ leaf of $T_{2^k}$, and the star $T_{2^l}$ to the other
leaves. Then contract all internal edges of the resulting tree, except the
internal edge of $T$. We have thus obtained $2^k$ trees $T_i$, $i=1,\ldots,2^k$,
with $2^{k+l}$ leaves and one internal edge, which define generators
$\alpha_{T_i}$ of $J_{2^{k+l}}$. Finally the second morphism is defined by
sending $\alpha_T$ to the product $\alpha_{T_1}\ldots\alpha_{T_{2^{k+l}}}$
(the factors commute among themselves). We omit the proof that this induces a
well-defined homomorphism.\\
 Moreover, it can be checked that
$\exp^{l}(\alpha_S)$ (where $\alpha_S$ is a generator of $J_{2^k}$) commutes
with $\alpha_{T_1}\ldots\alpha_{T_{2^{k+l}}}$, so that the two morphisms
commute. Indeed, notice that $\exp^l(\alpha_S)=\alpha_{\exp^l(S)}
\alpha_{(1,2^l)}\alpha_{(2^l+1,2^{l+1})}\ldots \alpha_{(2^{l+k-1}+1,2^{l+k})}$,
  where for simplicity we have assumed that $S$ corresponds to the labels
  $(1,\ldots,r)$, with $r\leq 2^k$, so that $\exp^l(S)$ corresponds to the
  labels $(1,\ldots,r2^l)$. Then
$$\exp^l(\alpha_S)\alpha_{T_1}= \alpha_{\exp^l(S)}
\alpha_{(1,2^l)}\alpha_{(2^l+1,2^{l+1})}\ldots
  \alpha_{(2^{l+r-1}+1,2^{l+r})}\alpha_{T_1}$$
$$=\alpha_{\exp^l(S)}\alpha_{(1,2^l)}\alpha_{T_1}\alpha_{(2^l+1,2^{l+1})}\ldots
  \alpha_{(2^{l+r-1}+1,2^{l+r})},$$
but
$$\alpha_{\exp^l(S)}\alpha_{(1,2^l)}\alpha_{T_1}=\alpha_{\exp^l(S)}\alpha_{j_{(1,2^l)}{T_1}}\alpha_{(1,2^l)}=\alpha_{T_r}\alpha_{\exp^l(S)}\alpha_{(1,2^l)}.$$
Since $\alpha_{T_i}$'s commute among themselves, we finally get 
$$\exp^l(\alpha_S)\alpha_{T_1}\ldots\alpha_{T_{2^{k+l}}}=\alpha_{T_1}\ldots\alpha_{T_{2^{k+l}}}\exp^l(\alpha_S).\;\;\square $$

Let now $M_{J}:=\coprod_{k\in \N} BJ_{2^k}$ be the associated topological monoid,
and, $(M_{J})_{\infty}\cong \Z\times
BJ_{2^{\infty}}$ be the inductive limit of the telescope induced by the right
multiplication $J_{2^k}\rightarrow J_{2^{k+1}}$ by $1\in J_2$.\\   

As for the spaces $\widetilde{\cal M}_{0,n+1}$, they are suitable models for the
classifying spaces $B(PJ_{n})$, and from the operadic structure maps
$$(*)\;\;\widetilde{\cal M}_{0,r+1}\times\widetilde{\cal M}_{0,n_1+1}\times\ldots\times
\widetilde{\cal M}_{0,n_r+1}\rightarrow \widetilde{\cal M}_{0, n_1+\ldots+n_r+1}$$ we get the
composition laws $\widetilde{\cal M}_{0,2^k +1}\times\widetilde{\cal M}_{0,2^l +1}\rightarrow
\widetilde{\cal M}_{0,2^{k+l}+1}$ and form the topological monoid
$M_{PJ}:=\coprod_{k\in\N} \widetilde{\cal M}_{0,2^k +1}$. In the
associated telescope, right multiplication by the point $\widetilde{\cal M}_{0,3}$ is
the expansion map $\widetilde{\cal M}_{0,2^k +1}\rightarrow \widetilde{\cal M}_{0,2^{k+1}
  +1}$. Similarly to the previous cases, we have $(M_{PJ})_{\infty}\cong
\Z\times B(PJ_{2^{\infty}})=\Z\times\widetilde{\cal M}_{0, \infty}$.\\

\noindent{\bf Consequences:} According to the ``Group-Completion" theorem of Quillen, the canonical map $(M_{?})_{\infty}\rightarrow
 \Omega B M_{?}$ of $H$-spaces (where $\Omega X$ denotes the based loop space
 of a pointed topological space $X$, and ? must be replaced by $PJ$,
 $J$ or $\Sigma$) is a strong homology equivalence (i.e. induces
 homology isomorphisms with any local system of abelian coefficients). In
 particular, the commutator subgroup of $\pi_1((M_{?})_{\infty})$ is
 perfect. To conclude, we have the

\begin{co}[Delooping $BPJ_{2^{\infty}}\cong\widetilde{\cal M}_{0, \infty}$ and $BJ_{2^{\infty}}$] 
There are homological equivalences $BPJ_{2^{\infty}}\rightarrow (\Omega
BM_{PJ})^0$, $BJ_{2^{\infty}}\rightarrow (\Omega
BM_{J})^0$ (where $(.)^0$ denotes the functor taking the connected component
of the unit of a monoid), and the commutator subgroups
$[PJ_{2^{\infty}},PJ_{2^{\infty}}]$ and $[J_{2^{\infty}},J_{2^{\infty}}]$ are
perfect.
\end{co}

\section{Concluding remarks}
 
The result of the last section may surely be improved: it must be considered
as a first step in delooping the towers and the quasi-braid groups. A related
question would be to detect the algebras over the graded operad
$\{H_*(PJ_n,\Z),\,n\geq 0\}$, having in mind that when considering  the
classical pure braid groups $P_n$, the algebras over the operad
$\{H_*(P_n,\Z),\,n\geq 0\}$ are the Gerstenhaber algebras, of great interest in physical mathematics.\\
Concerning the rest of this paper, many problems are open. The central one
concerns the relative natures of the three groups concerned: the
diffeomorphism group of the circle $Diff^+(S^1)$, Thompson group $T$, and
Neretin Spheromorphism group $N$, whose cohomologies resemble each other (the
continuous cohomology is considered for $Diff^+(S^1)$, and the $\Z/2\Z$
coefficients are natural for the group $N$ -- since anyway it is $\Q$-acyclic,
cf. \cite{kap2}). This can not be a
coincidence, and we would like to find a unified way to understand this
triangle of groups.

\end{document}